\input eplain
\input epsf.tex

 \input rotate.tex

\parskip \baselineskip
\parindent 0cm
\topmargin=1.2in
\bottommargin=1in
\leftmargin=1.25in
\rightmargin=1.25in

\newcount\notitre 
\newcount\nosection
\newcount \noproposition
\newcount \nolemma
\newcount \nocorollary

\font \bigtitlefont = cmbx12 at 14pt
\font \titlefont = cmbx12
\font \sectionfont = cmbx10 at 11pt 
\font \subsectionfont = cmbx10 at 11pt

\font \abstracttitle = cmbx10 at 11pt
\font \abstract = cmr10
\font \rm = cmr10 at 11pt
\font \sc = cmr8

\font \smallsf = cmss8
\font \tt = cmtt10
\font \smalltt = cmtt8

\def \titre  #1\par%
{ \bigskip\penalty -2000\nosection=0 \advance\notitre by 1
{\parindent=0 cm \titlefont \the\notitre. #1}\bigskip}

\def \section #1\par%
{\bigskip\advance\nosection by 1\bigskip%
{\parindent=0 cm \sectionfont \the\notitre.\the\nosection\ #1}\par}

\def \subsection #1\par{\medskip { \subsectionfont #1}\par}

\def \proposition #1\par{
\advance \noproposition by 1
\bigskip{\bf Proposition \the\notitre.\the\noproposition \ \ (#1)}\par
}

\def \lemma #1\par{
\advance \nolemma by 1
\smallskip{\bf Lemma \the\notitre.\the\nolemma \ \ (#1)}\par
}
\def \proof{\smallskip{\bf proof}\ }

\def \corollary #1\par{
\advance \nocorollary by 1
\smallskip{\bf Corollary \the\notitre.\the\nocorollary \ \ (#1)}\par
}

\font \mathbb=msbm10
\font \mathbbs=msbm8
\font \mathbbss=msbm6

\footline{\abstract \ifnum\pageno=1\else\hfil\the\pageno\hfil\fi}
\rm


\def \smallZ{\hbox{\mathbbs Z}}
\def \smallR{\hbox{\mathbbs R}}
\def \ssmallR{\hbox{\mathbbss R}}

\def\sup{\mathop{\hbox{\rm sup}}}
\def\exp{\mathop{\hbox{\rm exp}}}

\def \det{\mathop{\hbox{\rm det}}}

\def \br{\hfill\break }
\def \Np (#1){\|#1\|_p}
\def \Ndeux (#1){\|#1\|_2}

\def \ie{{\it i.e.\ }}
\def \R{\hbox{\mathbb R}}
\def \Z{\hbox{\mathbb Z}}
\def \N{\hbox{\mathbb N}}
\def \echn(#1){#1_1\ldots,#1_n}
\def \echd(#1){{#1^1\!\!\ldots,#1^d}}

\def \prodalpha {\alpha_{jk^1}\ldots\alpha_{jk^d}}

\def \pralpha #1,#2{\alpha_{j#1^1}^#2\ldots\alpha_{j#1^d}^#2}
\def \hatpralpha #1,#2{\hat\alpha_{j#1^1}^#2\ldots\hat\alpha_{j#1^d}^#2}

\def \hatprodalpha {\hat\alpha_{jk^1}\ldots\hat\alpha_{jk^d}}
\def \hatmulalpha {\hat\mulalpha}

\def \mulalpha {\alpha_{j\echd(k)}}

\def \Pnfa {P^n_{f_A}}
\def \Pnf {P^n_{f}}

\def \Enfa {E^n_{f_A}}

\def \Pfa {P_{f_A}}

\def \puncrem{\vskip 2pt\vbox{\hbox{\vrule height .8mm width .8mm}}\bigskip}
\def \endproof{
\vbox{\hrule\hbox to 5pt{\vrule height 4,2 pt\hfil\vrule}\hrule
}
\bigskip}

\belowlistskipamount 12pt
\interitemskipamount 10pt

\vglue .5in
\centerline{\bigtitlefont Independent Component Analysis by Wavelets}
\bigskip 
\bigskip

\leftline{{\bf Pascal Barbedor}} 
{\it
\leftline{Laboratoire de probabilit\'es et mod\'eles al\'eatoires}
\leftline{CNRS-UMR 7599, Universit\'e Paris VI et 
Universit\'e Paris VII}
\leftline{PO box 7012, 75251 PARIS Cedex 05 (FRANCE)}
\leftline{barbedor@proba.jussieu.fr}

} \bigskip
 
\bigskip
 
 \centerline{\abstracttitle Abstract}
{\abstract \advance \leftskip by .25in
\advance \rightskip by .25in 
We propose an ICA contrast based on the 
density estimation of the observed signal and its marginals 
by means of wavelets. 
The risk of the associated moment estimator is linked with 
approximation properties in Besov spaces. 
It is shown to converge faster than the at least 
expected minimax rate carried over from the underlying density estimations. 
Numerical simulations performed on some common types of densities
yield very competitive results, with a  high sensitivity to small departures from independence.

\br
 \par
}
\bigskip
 
\titre Introduction

In
signal processing, blind source separation consists in the
identification
of analogical, independent signals mixed by a black-box device.
In psychometry, one has the notion of structural latent variable 
whose mixed effects are only measurable through series of tests; an example 
are the Big Five (components of personality)
identified from factorial analysis 
by researchers in the domain of personality evaluation (\nocite{roch}Roch, 1995). 
Other application fields such as digital imaging, 
biomedicine, finance and econometrics also use  models aiming to 
recover hidden independent factors from observation. 
Independent component analysis (ICA) is one such tool; it can be seen as an extension
of principal component analysis, in that it goes beyond a simple 
linear decorrelation only satisfactory for a normal distribution; 
or as a complement, since its application is precisely pointless under the assumption of normality.

Papers on ICA are found in 
the research  fields of signal processing,  neural networks, 
statistics and information theory. Comon (1994)\nocite{ica-comon}
defined the concept of ICA as maximizing the degree of statistical independence among 
 outputs using contrast functions approximated by the Edgeworth expansion of 
 the Kullback-Leibler divergence.

The model is usually stated as follows: 
let  $x$ be a random variable on $\R^d$, $d \geq2$;   
one tries to find couples $(A,s)$, such that $x=As$, 
where $A$ is a square invertible matrix and $s$ 
a latent random variable whose components are mutually independent.
This is usually done through some contrast function that cancels out if and only if
the components of $Wx$ are independent, where $W$ is a candidate for the 
inversion of $A$.


Maximum-likelihood methods and contrast functions based on
mutual information or other divergence measures between densities are
commonly employed.
Cardoso (1999)\nocite{cardoso11} used higher-order cumulant tensors, which led to the Jade algorithm,
Bell and Snejowski (1990s) published an approach based on the Infomax principle. 
Hyvarinen and Oja (1997) presented the fast ICA algorithm\nocite{intro-hyva}\nocite{ICA-oja}.
\nocite{bell-snejI}

In the semi-parametric case, where the latent variable density
is left unspecified, 
Bach\nocite{kernel-ica} and Jordan (2002) proposed a contrast function based on canonical correlations in
a reproducing kernel hilbert space.
 Similarly, Gretton et al (2003)
\nocite{kmi} proposed kernel covariance and kernel mutual information
contrast functions.

The density model assumes that the observed random variable 
$X$ has the density $f_A$
given by 
$$
\eqalign{
f_A(x)&=|\det A^{-1}| f(A^{-1}x)\cr
&= |\det B| f^1(b_1x)\ldots f^d(b_dx),\cr
}$$
where $b_\ell$ is the $\ell^{th}$ row  of the matrix $B=A^{-1}$; 
this resulting from a change of variable if the latent density
$f$ is equal to the product of its marginals $f^1\ldots f^d$. 
In this regard, 
 latent variable $s=(s^1\!,\ldots,s^d)$ having independent components means the indepence of 
the random variables $s^\ell\circ \pi^\ell$ defined
on some product probability space $\Omega=\prod \Omega^\ell$, with $\pi^\ell$ the canonical projections.
So $s$ can be defined as the compound of the
unrelated $s^1,\!\ldots,s^d$ sources.  

Tsybakov and Samarov (2002) \nocite{tsybaica} proposed a method of simultaneous estimation 
of the directions $b_i$, based on nonparametric estimates of  
matrix functionals using the gradient of $f_A$.

In this paper, we propose a wavelet based ICA contrast.
The wavelet contrast $C_j$
compares the mixed density $f_A$ and its marginal distributions 
through their projections on a multiresolution analysis 
at level $j$. 
It thus relies upon the procedures of wavelet density estimation
which are found in a series of articles from Kerkyacharian and Picard (1992)\nocite{pikerk92} and
Donoho et al. (1996)\nocite{pic96}.

As will be shown,  the wavelet contrast has
the property to be zero only on a projected density with independent
components. The key parameter of the method lies in
the choice of a
resolution $j$, so that minimizing the contrast 
at that resolution yields a satisfactory approximate solution 
to the ICA problem.

The wavelet contrast can be seen as a special case of quadratic dependence measure, 
as presented in Achard et al. (2003)\nocite{achard2003}, which is equal to zero under 
independence. But in our case, the
resolution parameter $j$ allows more flexibility in controlling
the reverse implication.
Let's mention also that the idea of comparing in the $L_2$ norm a joint density
with the product of its marginals, can be traced
back to Rosenblatt (1975)\nocite{rosenblatt1975}.

Besov spaces are a general tool in describing smoothness properties
of functions; they also constitute the natural choice when dealing with
projections on a multiresolution analysis. 
We first show that a linear 
mixing operation is conservative as to Besov membership;
after what we are in position to derive a risk bound  
that will hold for the entire ICA minimization procedure.

Under its simplest form, the wavelet contrast estimator is a linear function
of the empirical measure on the observation. We give the rule for the
choice of a resolution level $j$  minimizing the risk, assuming
 a known regularity $s$ for a latent signal in some Besov space $B_{spq}$.

The estimator complexity is linear in the sample size
but exponential in the dimension $d$ of the problem; this is 
on account of an implicit  multivariate density estimation.
In compensation to this computational load,
the wavelet contrast shows a very good sensitivity 
to small departures from independence, and encapsulates
all practical tuning in a single parameter $j$.

\bigskip

\titre Notations

We set here the main notations and recall some definitions for  the convenience of ICA specialists. The reader already familiar with wavelets and Besov spaces can skip this part.

\unorderedlist

\li {\it Wavelets} 

Let $\varphi$ be some function of $L_2(\R)$ such that
the family of translates $\{\varphi(.-k),\, k\in \Z\}$ is an orthonormal 
system; let $V_j\subset L_2(\R)$ be the subspace spanned by 
 $\{\varphi_{jk}= 2^{j/2}\varphi(2^j . -k), k \in \Z\}$.

By definition, the sequence of spaces $(V_j), j\in \Z$,
is called a multiresolution analysis (MRA) of $L_2(\R)$ if
$V_j \subset V_{j+1}$ and $\bigcup_{j\geq 0} V_j$ is dense in $L_2(\R)$;
$\varphi$ is called the father wavelet or scaling function.

Let $(V_j)_{j\in \smallZ}$ be a multiresolution analysis
of $L_2(\R)$, with $V_j$ spanned by 
 $\{\varphi_{jk}= 2^{j/2}\varphi(2^j . -k), k \in \Z\}$. 
Define $W_j$  as the complement of  $V_j$ in $V_{j+1}$, and let the families
$\{\psi_{jk}, k\in \Z\}$ be a basis for $W_j$,
with $\psi_{jk}(x)= 2^{j/2}\psi(2^j x -k)$. 
Let  $\alpha_{jk}(f)=<f,\varphi_{jk}>$ and $\beta_{jk}(f)=<f,\psi_{jk}>$.
 
A function $f\in L_2(\R)$ admits a wavelet expansion on $(V_j)_{j\in \smallZ}$ if
 the series
$$\sum\limits_{k} \alpha_{j_0k}(f) \varphi_{jk}
+ \sum\limits_{j=j_0}^\infty
\sum\limits_{k} \beta_{jk}(f) \psi_{jk}
$$ is convergent to $f$ in $L_2(\R)$;
$\psi$ is called a mother wavelet.
\bigskip

The definition of a multiresolution analysis on $L_2(\R^d)$
follows the same pattern. 
But an MRA in dimension one also induces an associated 
MRA in dimension $d$, using the tensorial product procedure below.

Define $V^d_j$ as the tensorial product
of $d$ copies of $V_j$.
The increasing 
sequence $(V^d_j)_{j\in \smallZ}$
defines a multiresolution analysis of $L_2(\R^d)$
(Meyer, 1997)\nocite{meyer}:

for $(i^1\ldots,i^d) \in \{0,1\}^d$
and $(i^1\ldots,i^d) \neq (0\ldots,0)$,  
define $\Psi(x)_{i^1\ldots,i^d}=\prod_{\ell=1}^d \psi^{(i^\ell)}(x^\ell)$,
 with $\psi^{(0)}=\varphi$, $\psi^{(1)}=\psi$,
so that $\psi$ appears at least once in the product $\Psi(x)$
(we now on omit $i^1\ldots,i^d$ in the notation for $\Psi$, and
in \eqref{wavexpansion}, although it
is present each time);

for $(i^1\ldots,i^d) = (0\ldots,0)$,  define 
  $\Phi(x)=\prod_{\ell=1}^d \varphi(x^\ell)$;

for $j \in \Z$, $k\in \Z^d$, $x \in \R^d$, let  
$\Psi_{jk}(x)=2^{jd \over 2} \Psi(2^jx-k)$ and 
$\Phi_{jk}(x)=2^{jd \over 2} \Phi(2^jx-k)$;

define $W^d_j$ as the orthogonal complement of  
$V^d_j$ in $V^d_{j+1}$; 
it is an orthogonal sum of $2^d-1$ spaces
having the form $U_{1j}\ldots\otimes U_{dj}$, where $U$
is a placeholder for $V$ or $W$; $V$ or $W$ 
are thus placed using up all permutations, but with 
  $W$ represented at least once, so that 
a fraction of the overall innovation brought by
 the finer resolution $j+1$ is always present in the tensorial product.

A function $f$ admits a wavelet expansion on the basis $(\Phi, \Psi)$ 
if the series
$$\sum\limits_{k\in \smallZ^d} \alpha_{j_0k}(f) \Phi_{j_0k}
+ \sum\limits_{j=j_0}^\infty \sum\limits_{k\in \smallZ^d} \beta_{jk}(f) 
\Psi_{jk},
\eqdef{wavexpansion}$$ is convergent to $f$ in $L_2(\R^d)$.

In fact, with the concentration condition 
$$\sum_k |\varphi(x+k)| \leq C\; a.s.,\eqdef{concentration condition}$$
verified in particular for a compactly supported wavelet, 
any function in $L_1(\R^d)$ admits a wavelet expansion.
Otherwise any function in a Besov space $B_{spq}(\R^d)$ admits a wavelet expansion.

\bigskip

In connection with function approximation, wavelets can be viewed 
as falling in the category of orthogonal
series methods, or also in the category of kernel methods.

The approximation at level $j$ of a funtion $f$ that admits a 
multiresolution
expansion is the orthogonal projection  $P_jf$ of $f$ onto
$V_j\subset L_2(\R^d)$  defined by:
$$(P_jf)(x)=\sum_{k \in \smallZ^d}\alpha_{jk}\Phi_{jk}(x),
$$
where  
 $\alpha_{jk}=\mulalpha=\int f(x)\Phi_{jk}(x)\,dx$.

With the concentration condition above, 
the projection operator can also be  written
$$
(P_jf)(x)=\int_{\smallR^d} K_j(x,y)f(y)d(y),
$$
with $K_j(x,y)=2^{jd}\sum_{k\in \smallZ^d}\Phi_{jk}(x-k)
\overline{\Phi_{jk}(y-k)}$. $K_j$ is an orthogonal projection kernel
 with window
$2^{-jd}$ (which is not translation invariant).

\bigskip

\li {\it Besov spaces} 

Let $f$ be a function in $L_p(\R^d)$ and $h\in \R^d$.
Define the first order difference $\Delta_hf$ by $\Delta_hf(x)=f(x+h)-f(x)$
and the $k^{th}$ order difference $\Delta^k_hf=\Delta_h\Delta^{k-1}_hf$ ($k=1, 2,\ldots$ with  $\Delta_h^0f=f,\; \Delta_h^1f=\Delta_hf$).

The modulus of continuity of order $k$  of $f$ in the metric of $L_p$, along direction $h$,
is defined by (Nikol'ski$\breve{\hbox{\i}}$, 1975, p.145-160)
$$
\omega^k_h(f,\delta)_p=\sup_{|t|\leq \delta} \|\Delta^k_{th}f(x)\|_p,\qquad \delta \geq0,\quad |h|=1.
$$

The modulus of continuity of order $k$ of $f$ in the direction of 
the subspace $\R^m\subset \R^d$ is defined by
$$
\Omega^k_{\smallR^m}(f,\delta)_p=\sup_{|h|=1, h\in\smallR^m}
\omega^k_h(f,\delta)_p.
$$

If the function $f$ has arbitrary derivatives of order $\varrho$ relative to the first $m$ coordinates,
one can define, for $h\in \R^m$,$$
f^{(\varrho)}_h=\sum_{|n|=\varrho} f^{(n)}h^n,
$$ with $h=(h_1,\ldots,h_m,0,\ldots, 0)$, $|h|=1$, $|n|=\sum_1^mn_i$ and
$h^n=h_1^{n_1}\ldots h_m^{n_m}=h_1^{n_1}\ldots h_m^{n_m}0^0\ldots0^0$.

The modulus of continuity of order $k$ of the derivatives of order $\varrho$
 of $f$ is then defined by
$$\Omega^k_{\smallR^m}(f^{(\varrho)},\delta)_p= 
\sup_{|h|=1, h\in\smallR^m}
\omega^k_h(f^{(\varrho)}_h,\delta)_p=\sum_{|n|=\varrho} \Omega^k_{\smallR^m}(f^{(n)},\delta)_p.$$

Let $s=[s]+\alpha$; the H\" older space $H^s_p(\R^d)$ is defined as the collection
of functions in $L_p(\R^d)$ such that
$$\eqalign{
\|\Delta_hf^{(n)}\|_p \leq M|h|^\alpha,
\quad
\forall n=(n^1,\ldots,n^d),\; \hbox{ with } |n|=\sum_1^dn_i=[s],\cr
\hbox{or equivalently, }
\Omega_{\smallR^d}(f^{([s])},\delta)_p=\sup_{h\in \smallR^d}\omega_h(f^{([s])},\delta)_p \leq M\delta^\alpha,
}$$ where $M$ does not depend on $h$.

Besov spaces introduce  a finer scale  of smoothness than 
is provided by H\" older spaces. For each $\alpha>0$
this can be accomplished  by introducing a second parameter $q$ and
applying ($\alpha$, $q$) quasi-norms  (rather than ($\alpha$, $\infty$))
to the modulus of continuity of order $k$.

Let $s>0$ and $(\varrho,k)$ forming an admissible  pair
of nonnegative integers satisfying the inequalities 
$k>s-\varrho>0$.
By definition, $f \in L_p(\R^d)$ belongs to the class $B_{spq}(\R^d)$
if there exist generalized partial derivatives of
$f$ of order $n=(n^1,\ldots,n^d)$, $|n|\leq \varrho$,
and  one of the following semi-norms is finite:
$$\eqalign{
J'_{spq}(f)=\sum_{|n|=\varrho}\left(  
\int_0^\infty | t^{-(s-\varrho)} \Omega^k_{\smallR^d}(f^{(n)},t)_p|^q {dt\over t}
\right)^{1\over q}, \cr
J''_{spq}(f)=\left(  
\int_0^\infty | t^{-(s-\varrho)} \Omega^k_{\smallR^d}(f^{(\varrho)},t)_p|^q {dt\over t}
\right)^{1\over q}.}\eqdef{Jspqnikol}
$$

For fixed $s$ and $p$, the space $B_{spq}$ gets larger with increasing $q$.
In particular, for $q=\infty$, $B_{spq}(\R)=H^s_p(\R)$;
various other 
embeddings 
exist since Besov spaces cover many well known classical concrete
function spaces having their own history.

Finally, Besov spaces also admit a characterization in terms of wavelet coefficients,
which makes them intrinsically connected
to the analysis of curves via wavelet techniques. 

\filbreak
 $f$ belongs to the (inhomogeneous) Besov space $B_{spq}(\R^d)$ if
$$J_{spq}(f)=\| \alpha_{0.}\|_{\ell_p} + 
\left[ \sum_{j\geq 0}\left[  
2^{js}2^{dj({1\over 2}-{1\over p}})\|\beta_{j.} \|_{\ell_p}
\right]^q\right]^{1\over q} < \infty,\eqdef{jspq}$$
with  $s>0$, $1\leq p \leq\infty$, $1\leq q \leq\infty$, 
and $\varphi, \psi$ $\in C^r, r>s$ (Meyer, 1997).

\endunorderedlist

A more complete presentation of wavelets linked with Sobolev and Besov approximation theorems 
and statistical applications can be found in 
the book from H\"ardle et al. (1998)\nocite{picardbook}. 
General references about Besov spaces
are  Peetre (1975), 
Bergh \& L\"ofstr\"om (1976), Triebel (1992), DeVore \& Lorentz (1993). 
\bigskip

\titre Wavelet contrast, Besov membership 
 
\def \fix(#1,#2){#1^{\star#2}}

Let $f$ be a density function with marginal 
 distribution in dimension $\ell$,
$$ x^\ell \mapsto \int_{\smallR^{d-1}}
f(\echd(x))\,dx^1\ldots dx^{\ell-1}dx^{\ell+1}\ldots dx^d,$$
denoted by $\fix(f,\ell)$.

As integrable positive functions, $f$ and the $\fix(f,\ell)$
 admit a wavelet expansion 
on a  basis $(\varphi,\psi)$ verifying the concentration condition \eqref{concentration condition}.
One can then consider the projections up to order $j$,  
that is to say the projections of $f$ and $\fix(f,\ell)$ 
on  $V^{d}_j$ and $V_j$ respectively, namely
$$
P_jf(x)=\sum_{k \in \smallZ^d}\alpha_{jk}\Phi_{jk}(x)
 \ \ \hbox{and}\ \ 
P_j^\ell\fix(f,\ell)(x^\ell)=\sum_{k^\ell \in \smallZ}\alpha_{jk^\ell}\varphi_{jk^\ell}(x^\ell),
$$
where  
$\alpha_{jk^\ell}=\int \fix(f,\ell)(x^\ell)\varphi_{jk^\ell}(x^\ell)\,dx^\ell$
and $\alpha_{jk}=\mulalpha=\int f(x)\Phi_{jk}(x)\,dx$.

 
\proposition wavelet contrast\par
\definexref{wavelet contrast}{\the\notitre.\the\noproposition}{proposition}
{
\it
\noindent
Let $f$ be a density function on $\R^d$ and let $\varphi$ be the
scaling function of a multiresolution analysis
verifying the concentration condition \eqref{concentration condition}.

Define the contrast function 
$$C_j(f)=\sum\limits_{\echd(k)}(\mulalpha - \prodalpha)^2,$$
with  
$\alpha_{jk^\ell}=\int_{\ssmallR} \fix(f,\ell)(x^\ell)\varphi_{jk^\ell}(x^\ell)\,dx^\ell$
and $\mulalpha=\int_{\ssmallR^d} f(x)\Phi_{jk^1,\ldots,k^d}(x)\,dx$.

The following relation hold: 
$$
 f \hbox{ factorizable} \kern 7pt\Longrightarrow \kern 5pt C_j(f)=0.$$

If $f$ and $\varphi$ are compactly supported or else if $f\in L_2(\R^d)$,
the following relation hold:
$$
  C_j(f)=0 \kern 7pt\Longrightarrow \kern 5pt P_jf = \prod_{\ell=1}^d P^\ell_j\fix(f,\ell).
$$

}

{\bf Proof}

As for the first assertion, with $f = f^1\ldots f^d$, one has  $\fix(f,\ell)=f^\ell$, $\ell=1,\ldots d$. 
Whence for $k=(k^1,\ldots, k^d)\in \Z^d$, one has by Fubini theorem, 
$$\eqalign{\alpha_{jk}(f)&=\alpha_{jk}(\fix(f,1)\ldots \fix(f,d))
=\int_{\R^d}\fix(f,1)\ldots \fix(f,d) \Phi_{jk}(x) dx\cr
&= \int_{\R}\fix(f,1)\varphi_{jk^1}(x^1)dx^1 \ldots\int_{\R}\fix(f,d)\varphi_{jk^d}(x^d)dx 
=\alpha_{jk^1}(\fix(f,1))\ldots\alpha_{jk^d}(\fix(f,d)).
}
$$

For the second assertion,  $C_j = 0$ entails
$\alpha_{jk}(f)=\alpha_{jk^1}(\fix(f,1))\ldots\alpha_{jk^d}(\fix(f,d))$ for all $k\in \Z^d$.
So that for $P_jf\in L_p(\R^d)$,
$$\eqalign{P_jf&=\sum_k\alpha_{jk}(f)\Phi_{jk}= \sum_k 
\alpha_{jk^1}(\fix(f,1))\varphi_{jk^1}\ldots\alpha_{jk^d}(\fix(f,d))\varphi_{jk^d}\cr
&=\sum_{k^1}\alpha_{jk^1}(\fix(f,1))\varphi_{jk^1}\ldots 
\sum_{k^d}\alpha_{jk^d}(\fix(f,d))\varphi_{jk^d}\cr
&=P_j^1\fix(f,1)\ldots P_j^d\fix(f,d),
}
$$ with passage to line 2 justified by the fact that
$(\alpha_{jk}(f)\Phi_{jk})_{k\in \smallZ^d}$ 
is a summable family of $L_2(\R^d)$ or else  is a finite sum for
a compactly supported density and a compactly supported wavelet.

\endproof


For the zero wavelet contrast to give any clue as to whether     
the non projected difference $f-\fix(f,1)\ldots \fix(f,d)$ is itself close to zero, a key parameter lies in the order of projection $j$.

Under the notations of the preceding proposition, with a zero wavelet contrast and  
assuming existence in $L_p$, one has $\| P_{j}f- P^1_{j}\fix(f,1) \ldots P^d_{j}\fix(f,d)\|_p=0$, 
and so
$$\eqalign{
\|\, f- \fix(f,1)\ldots\fix(f,d) \,\|_p
&\leq \| f-   P_j f\|_p
+ \| P^1_{j}\fix(f,1)\ldots P^d_{j}\fix(f,d) - \fix(f,1)\ldots\fix(f,d)\|_p\cr
&= \| f-   P_j f\|_p
+ \| P_{j}(\fix(f,1)\ldots \fix(f,d)) - \fix(f,1)\ldots\fix(f,d)\|_p.\cr
}$$

If we now impose some regularity conditions on the densities, 
in our case if we now require that  $f$ and the product
of its marginals belong to the (inhomogeneous) Besov space $B_{spq}(\R^d)$,
 the approximation error can be 
evaluated precisely.
With a $r$-regular wavelet $\varphi$, $r > s$, the very definition 
of Besov spaces implies for any member $f$ that (Meyer, 1997)
$$
\| f-   P_j f\|_p = 2^{-js}\, \epsilon_j, \quad \{\epsilon_j\} \in \ell_q(\N^d).\eqdef{meyer}
$$

{\it Remark}

In the special case where $f_A$ and the product of its marginals belong to $L_2(\R^d)$, Parseval equality implies  
that $C_j$ is equal to the square of the $L_2$ norm of 
$P_jf_A -P_j^1\fix(f_A,1)\ldots P_j^d\fix(f_A,d)$.
And one can write, 
$$
\eqalign{C_j(f_A)^{1\over2}  &= \| P_{j}(\fix(f_A,1)\ldots \fix(f_A,d))-   P_j f_A\|_2\cr
&\leq \| f_A-   P_j f_A\|_2 
+  \| f_A- \fix(f_A,1)\ldots \fix(f_A,d)\|_2 
+ \| P_{j}(\fix(f_A,1)\ldots \fix(f_A,d)) - \fix(f_A,1)\ldots\fix(f_A,d)\|_2,\cr}
$$
hence with notation $K_\star(A,f) = \| f_A- \fix(f_A,1)\ldots \fix(f_A,d)\|_2$,
$$  |K_\star(A,f)- C_j(f_A)^{1\over2}| \leq 
 \| f_A-   P_j f_A\|_2 
+ \| P_{j}(\fix(f_A,1)\ldots \fix(f_A,d)) - \fix(f_A,1)\ldots\fix(f_A,d)\|_2,\eqdef{ketoile}
$$
which gives another illustration of the shrinking with $j$ distance between 
the wavelet contrast and the true norm evaluated at $f_A$.
In particular when $A\neq I$, $C_j(f_A)$ cannot be small and for
$A=I$, $C_j$ must be small, for $j$ big enough.

Continuing on the special case $p=2$, the wavelet contrast can be viewed as an example of quadratic 
dependence measure
as presented in the paper from Achard et al (2003).

Using the orthogonal projection kernel associated to the function $\varphi$,
 one has the writing
$$C_j(f_A)=\int_{\R^d}\left(\Enfa K_j(x,Y)-\prod_{i=1}^d \Enfa K^i_j(x^i,Y^i)\right)^2dx,$$
with $K_j(x,y)=2^{jd}\sum_{k\in\smallZ^d}\Phi_{jk}(x-k)\Phi_{jk}(y-k)$
and $K^i_j(x,y)=2^j\sum_{k\in\smallZ}\varphi_{jk}(x^i-k^i)\varphi_{jk}(y^i-k^i)$.

This is the form of the contrast in the paper from Achard et al. (2003),
except that in our case the kernel is not scale invariant;
but the ICA context is scale invariant by feature, since
the inverse of $A$ is conventionally determined up to a scaling diagonal or permutation matrix, after a whitening step.
\puncrem

To take advantage of  relation \eqref{meyer} in the ICA context, we need a fixed Besov space containing the mixed density $f_A$ and the product of its marginals, for any invertible matrix $A$.

The two following propositions check that the mixing by $A$ is
conservative as to Besov membership, and that
the product of the marginals of a density $f$
belongs to the same Besov space than $f$. 
It is equivalent to assume that $f$ is in $B_{spq}(\R^d)$
or that the factors $f^i$ are in $B_{spq}(\R)$. If the factors 
have different Besov parameters, 
one can theoretically always find a bigger Besov space using Sobolev
inclusions
$$\eqalign{
B_{s'pq'} \subset B_{spq} &\quad\hbox{ for } s' \geq s,\quad q' \leq q; \cr
B_{spq} \subset B_{s'p'q} &\quad \hbox{ for }
p \leq p' \hbox{ and } s' = s + d/p' -d/p. 
}$$


\proposition Besov membership of marginal distributions\par
\definexref{marginal distributions membership}{\the\notitre.\the\noproposition}{proposition}
{\it
If $f$ is a density function belonging to $B_{spq}(\R^d)$ then
each of its marginals  belong to $B_{spq}(\R)$.
}

\proof

Let us first check the $L_p$ membership of the marginal distribution.
For $p\geq 1$, by convexity one has,  
$$
\int_{\smallR^d} | f_A |^p \,dx = \int_{\smallR} \int_{\smallR^{d-1}} |f_A |^p \,\fix(dx,\ell) dx^\ell
 \geq  \int_{\smallR} \left\vert \int_{\smallR^{d-1}} f_A \,\fix(dx,\ell) \right\vert^p dx^\ell
=  \int_{\smallR} | \fix(f_A, l)|^p \,dx^\ell;
$$that is to say  $\| \fix(f_A, l)\|_p \leq \|f_A\|_p$.

With  
the $\ell^{th}$ canonical vector of $\R^d$ denoted by $e^\ell$ and for $t\in \R$,
one has, 
$$\eqalign{
\Delta^k_t\fix(f,\ell)(x^\ell)
&= \sum_{l=0}^k (-1)^{l+k}C_k^l \fix(f,\ell)(x+t)
=\sum_{l=0}^k (-1)^{l+k}C_k^l \int_{\smallR^{d-1}} f(x+te^\ell) dx^{\ast\ell}
= \int_{\smallR^{d-1}} \Delta^k_{te^\ell}f(x) dx^{\ast\ell};
}
$$
so that
$$\eqalign{
\|\Delta^k_t\fix(f,\ell)\|_{L_p(\smallR)}^p
&=\int_{\smallR} \left|\int_{\smallR^{d-1}} \Delta^k_{te^\ell}f(x)
 dx^{\ast\ell}\right|^pdx^\ell
\leq
 \int_{\smallR^d}\left| 
 \Delta^k_{te^\ell}f(x)\right|^p dx
 \leq \|\Delta^k_{te^\ell}f\|_{L_p(\smallR^d)}^p,  
}
$$
and  
$$\eqalign{
\omega^k(\fix(f,\ell),\delta)_p  
&= \sup_{|t| \leq \delta}
 \|\Delta^k_{t}\fix(f,\ell)\|_{Lp(\smallR)}
\leq 
\sup_{|t| \leq \delta} \|\Delta^k_{t e^\ell} f\|_{L_p(\smallR^d)}
 = \omega^k_{e^\ell}(f,\delta)_p,
}
$$
and
$$
\Omega^k(\fix(f,\ell),\delta)_p=
\omega^k(\fix(f,\ell),\delta)_p
\leq \omega^k_{e^\ell}(f,\delta)_p
\leq \sup_{|h|=1, h\in\smallR^d}
\omega^k_h(f,\delta)_p=\Omega^k_{\smallR^d}(f,\delta)_p.
$$

Using the admissible pair $(k,\varrho)=([s]+1,0)$, one
can see from \eqref{Jspqnikol} that $J'_{spq}(\fix(f,\ell))\leq J'_{spq}(f)$.

\endproof

Next, we check that the mixed density $f_A$ belongs to the
same Besov space than the original density $f$.

\proposition Besov membership of the mixed density\par
\definexref{mixed-membership}{\the\notitre.\the\noproposition}{proposition}

{\it
Let $f=f^1\ldots f^d$ and $f_A(x)=|\det A^{-1}| f(A^{-1}x)$.

(a) if $f \in L_p(\R^d)$, or if each
$f^\ell$ belongs to $L_p(\R)$, then $f_A$ and the product
 $\prod f^{\star\ell}_A$ belong to $L_p(\R^d)$.

(b) $f$ and $f_A$ have  same Besov semi-norm up to a constant.

Hence $f$ and $f_A$  belong to the same 
(inhomogeneous) Besov space $B_{spq}$.
}
\proof

For (a), with $p\geq 1$, as in Prop. 
\ref{marginal distributions membership} above,
one has  $\| \fix(f_A, \ell)\|_p \leq \|f_A\|_p$.
Also,
with  the determinant of $A$ denoted by $\vert A\vert$, 
$$\|f_A\|_p= \vert A\vert^{-p}\int \vert f(A^{-1}x)\vert^p \,dx 
= \vert A\vert^{-p}\int \vert f(x)\vert^p \,\vert A\vert\, dx
= \vert A\vert^{1-p}\,\|f\|_p.
$$
And finally by Fubini theorem, $\|f\|_{L_p(\smallR^d)}=\|f^1\|_{L_p(\smallR)} \ldots \|f^d\|_{L_p(\smallR)}$,
 so that $f\in L_p(\R^d) \Longleftrightarrow f^\ell \in L_p(\R), \ell=1\ldots d$ .

For (b), 
with  a change of variable in the integral one has, 
$$
\|\Delta_{th}f_A\|_p=
\vert A\vert ^{-1+ {1\over p}}\, \|\Delta_{tA^{-1}h} f\|_p\,;
$$ 
so that 
$$
\omega_h(f_A,\delta)_p=\sup _{|t|\leq \delta, |h|=1}
\|\Delta_{th}f_A\|_p
=|A|^{-1+{1\over p}}\sup _{|t|\leq \delta|A^{-1}h|, |h|=1}
\|\Delta_{th}f\|_p=\omega_l(f,\delta|A^{-1}h|)_p,
\quad|h|=1; 
$$
and $$
\Omega_{\smallR^d}(f_A,\delta)_p
=|A|^{-1+{1\over p}}\Omega_{\smallR^d}(f,\delta|A^{-1}h|)_p,\quad |h|=1.
$$

Next, with the change of variable $u=t|A^{-1}h|$, 
$$\eqalign{
\int_0^\infty |t^{-\alpha}\Omega(f_A,t)|^q {dt\over t}
&=(|A|^{-1+{1\over p}}\,|A^{-1}h|^\alpha)^q
\int_0^\infty |u^{-\alpha}\Omega(f,u)|^q {du\over u},
\quad |h|=1\cr
&\leq
(|A|^{-1+{1\over p}}\,\|A^{-1}\|^\alpha)^q
\int_0^\infty |u^{-\alpha}\Omega(f,u)|^q {du\over u}\,.
}%
$$

In view of \eqref{Jspqnikol}, using the 
admissible pair $(k,\varrho)=([s]+1,[s])$ yields
the desired result  when $0<s<1$.

When $1 \leq s$, with the same admissible pair $(k,\varrho)=([s]+1,[s])$, and by recurrence, since  
$df_A(h) = \vert A^{-1}\vert\, df(A^{-1}h)\circ A^{-1}$ one can see in the same way that
the modulus of continuity of
the (generalized) derivatives of $f_A$ or order $k$  are
bounded by those of $f$. 

Note that if $A$ is whitened, in the context of ICA, 
the norms of $f$ and $f_A$ are equal, at least when $s<1$.
\endproof

\titre Risk upper bound

Define the  experiment
${\cal E}^n=({\cal X}^{\otimes n},\, {\cal A}^{\otimes n},\,
(X_1, \ldots, X_n),\,
\Pnfa,\, f_A \in B_{spq})$, where $X_1, \ldots, X_n$ is an iid sample of $X=AS$, 
and $\Pnfa=\Pfa\ldots \otimes \Pfa$
is the joint distribution of ($\echn(X)$).
 Likewise, define $\Pnf$ as the joint distribution of $(\echn(S))$.

The coordinates $\alpha_{jk}$ in the wavelet contrast are estimated as usual by,
$$\hat \mulalpha = 
{1 \over n} \sum\limits_{i=1}^n \varphi_{jk^1}(X_i^1)\ldots\varphi_{jk^d}(X_i^d)
\hbox{ and }
\hat \alpha_{jk^\ell} = 
{1 \over n} \sum\limits_{i=1}^n \varphi_{jk^\ell}(X_i^\ell), \ \ell=1,\ldots d.
\eqdef{moments estimators}
$$
The linear wavelet contrast estimator is given by,
$$\hat C_j(\echn(x))=\sum\limits_{\echd(k)}
(\hatmulalpha - \hatprodalpha)^2
 = \sum\limits_{k\in \smallZ^d} \hat\delta_{jk}^2,\eqdef{deltajk}
$$
where we define $\hat\delta_{jk}$ as the difference $\hatmulalpha - \hatprodalpha$.

The estimator $\hat\alpha_{jk}$is unbiased under $\Enfa$, but so is
not $\hatprodalpha$ unless $A=I$, although it is asymptotically unbiased.

We also make the assumption that both the density and the wavelet
are compactly supported so that all sums in $k$ are finite. For
simplicity we further suppose the density support to be the hypercube, so that
$\sum_k 1 \approx 2^{jd}$.

To bound the risk we use a single lemma whose proof relies
on a classical U-statistic lemma, namely the connection between a U-statistic and its associated Von Mises statistic. To fit our purpose  the U-statistic lemma first needed to be adapted to kernels that are (generally unsymmetric) products of  
$\Phi_{jk}$ and $\varphi_{jk}\circ\pi^\ell$, and thus 
depend on the resolution parameter $j$. This is done in lemmas
\ref{UV} appearing  in the Appendix.

\filbreak
\lemma Moments of wavelet coefficients estimators\par
\definexref{various expectations}{\the\notitre.\the\nolemma}{lemma}
{\it
Let $\rho,\sigma\geq 0$ be fixed integers;
the following relations hold:
$$\eqalign{
\Enfa \hat\alpha_{jk}^\rho(\hatprodalpha)^\sigma= \alpha_{jk}^\rho(\prodalpha)^\sigma+O(n^{-1}).
}$$
And in corollary,
$\Enfa \hat\delta_{jk}^\rho=\delta_{jk}^\rho+O(n^{-1})$.

}%
\proof

To the statistic $V_{nj}=\hat\alpha_{jk}^\rho(\hatprodalpha)^\sigma$ corresponds a U-statistic $U_{nj}$ with unsymmetric kernel
 $$ h_j(x_1,\ldots,x_{\rho+d\sigma})= \Phi_{jk}(x_1)\ldots \Phi_{jk}(x_\rho)
 \varphi_{jk^{\ell_1}}\circ\pi^{\ell_1}(x_{\rho+1})
 \ldots\varphi_{jk^{\ell_{d\sigma}}}\circ\pi^{\ell_{d\sigma}}(x_{\rho+d\sigma}),$$
 with $\pi^\ell$ the canonical projection on component $\ell$,
 $(\ell_{i\sigma+1},\ldots,\ell_{(i+1)\sigma})=(1,\ldots,d)$, $i=0\ldots d-1$
 and $\Phi(x)=\prod_{\ell=1}^d\varphi\circ\pi^\ell(x)$.

By application of lemma \ref{UV} in Appendix,
$$|\Enfa \hat\alpha_{jk}^\rho(\hatprodalpha)^\sigma 
- \alpha^\rho(\prodalpha)^\sigma|
\leq \Enfa |V_{nj}-U_{nj}|=O(n^{-1})
$$
\endproof

We now express a risk bound for the wavelet contrast estimator. In particular
we show that the bias of the estimator is of 
the order of $C2^{jd}/ n$. This is better than the convergence
rate of $n^{-1\over2}$ for the empirical
Hilbert-Schmidt independence criterion
 (Gretton et al. 2004, theorem 3), except that 
 in our case 
 the resolution parameter $j$ must still be set to some value,
 especially 
 to cope with the antagonist objectives of  
 reducing the estimator bias and variance.

\proposition Risk upper bound for $\hat C_j$\par
\definexref{riskofCj}{\the\notitre.\the\noproposition}
{proposition}{\it 

The quadratic risk $\Enfa(\hat C_j -C_j)^2$ 
 and the bias $\Enfa \hat C_j - C_j$ have convergence rate $2^{jd}O(1/n)$.

In corollary, 
  the variance 
$\Enfa \left(\hat C_j - \Enfa\hat C_j\right)^2$ has convergence rate $ 2^{jd}O(1/n)$ and 
 the quadratic risk around zero is $ \Enfa \hat C_j^2 = C_j^2+ 2^{jd}O(1/n).$
}

\proof 

The risk about $C_j$ is written,
$$\Enfa \left( \hat C_j - C_j \right)^2
=\Enfa \sum_{k,\ell}
( \hat\delta_{jk}^2 - \delta_{jk}^2)
( \hat\delta_{j\ell}^2 - \delta_{j\ell}^2)
\eqdef{crossrisk}$$

For the squared terms where $k=\ell$,
lemma \ref{various expectations} yields directly
$\Enfa ( \hat\delta_{jk}^2 - \delta_{jk}^2)^2 = O(n^{-1})$,
so that the corresponding risk component is bounded
by $C2^{jd}n^{-1}$.

For crossed terms where $k\neq \ell$, observe that with a compactly supported Daubechies Wavelet $D2N$, whose support
is $[0,2N-1]$, only a thin band of terms around the diagonal is non zero: 
$$\varphi_{jk}\varphi_{j\ell}=0,\quad \hbox{ for }  |\ell-k|\geq 2N-1.$$

When  $k$ is fixed,
the cardinal of the set $|\ell -k|\leq 2N-1$ is lower than  $(4N)^d$; hence, by Cauchy-Schwarz inequality
and lemma \ref{various expectations},
$$\eqalign{
\Enfa \sum_{k,\ell}
(\hat\delta_{jk}^2 - \delta_{jk}^2)(\hat\delta_{j\ell}^2 - \delta_{j\ell}^2)
&\leq \sum_k\Enfa{}^{1\over2}(\hat\delta_{jk}^2 - \delta_{jk}^2)^2
\sum_{|\ell -k|\leq 2N-1}[
\Enfa(\hat\delta_{j\ell}^2 - \delta_{j\ell}^2)^2]^{1\over 2}\cr
&\leq 2^{jd} Cn^{-1/2}\;(4N)^d Cn^{-1/2}=C2^{jd}n^{-1}.
}$$

The bias convergence rate is a direct consequence
of lemma \ref{various expectations}, and the
two remaining assertions
follow from the usual relations, 
$\Enfa(\hat C_j -C_j)^2=\Enfa (\hat C_j - \Enfa\hat C_j)^2 + (\Enfa \hat C_j - C_j)^2$;
and $\Enfa \hat C_j^2=( \Enfa\hat C_j)^2 + \Enfa (\hat C_j - \Enfa\hat C_j)^2$.

\endproof

We now give a rule for choosing the resolution $j$ 
 minimizing the (about zero) risk  upper bound. 
This rule, obtained as usual through bias-variance balancing, 
depends on $s$, 
the unknown regularity of $f$, supposed to be a member of some Besov space $B_{spq}$.
The associated convergence rate improves upon the minimax  $n^{-2s\over 2s+d}$ of the 
underlying density estimations (see Kerkyacharian \& Picard,  1992).

\proposition minimizing resolution in the class $B_{s2\infty}$\par
\definexref{linear-choice}{\the\notitre.\the\noproposition}{proposition}
{\it
Assume that $f$ belongs to $B_{s2\infty}(\R^d)$,
and  that $C_j$ is based on a $r$-regular wavelet $\varphi$, $r > s$.

The minimizing resolution $j$ is such that $2^j\approx  n^{{1\over 4s + d}}$ and ensures a quadratic risk converging to zero at rate  
$n^{-{4s\over 4s +d}}$.
}

\proof

By Prop. \ref{marginal distributions membership} and 
\ref{mixed-membership} we know that $f_A$ and the product of its marginal
distributions belong to the same Besov space than the original $f$, so that equation \eqref{ketoile} becomes
$$ |K_\star(A,f)- C_j(f_A)^{1\over2}| \leq  K 2^{-js};\eqdef{KmoinsC}$$ 
with $K$ a constant.

Taking power 4 of \eqref{KmoinsC} and using prop. \ref{riskofCj},
$$\eqalign{
R(\hat C_j,f_A) 
 + K^\star Q(C_j^{1\over2}, K^\star) &\leq K 2^{-4js} + 2^{jd}Kn^{-1}\cr
},$$ with $K$ a placeholder for an unspecified constant, 
$Q(a,b) =-4a^3 + 6a^2b - 4ab^2 + b^3$, and $R$ denoting the quadratic risk around zero.

When $A$ is far from $I$, the constant $K_\star$ is strictly positive 
and the risk relative to zero has no useful upper bound.
Although the risk relative to $C_j$ is always 
in $2^{jd}Kn^{-1}$.

With $A$ getting closer to $I$, $K_\star$ is brought down to zero
and
the bound is minimum when, constants apart, we balance 
$ 2^{-4js}$ with $2^{jd}  n^{-1}$,
or $2^{j(d+4s)}$ with $n$.

This yields $2^j = O(n^{1\over 4s + d})$
and  convergence rate  $n^{-4s\over 4s+d}$ for the risk
relative to zero under independence and also
 for the risk relative to $C_j$ by substitution in the expression 
 given by Prop. \ref{riskofCj}.

\endproof

\corollary minimizing resolution in the class $B_{spq}$\par
\definexref{linear-choice}{\the\notitre.\the\noproposition}{proposition}
{\it
Assume that $f$ belongs to $B_{spq}(\R^d)$,
and  that $C_j$ is based on a $r$-regular wavelet $\varphi$, $r > s'$.

The minimizing resolution $j$ is such that $2^j\approx  n^{{1\over 4s' + d}}$, with $s'=s+d/2-d/p$ if $1\leq p\leq 2$ and $s'=s$
if $p>2$.

This resolution ensures a quadratic risk converging to zero at rate  
$n^{-{4s'\over 4s' +d}}$.

}

\proof

If $1 \leq p \leq 2$, using 
the Sobolev embedding $B_{spq} \subset B_{s'p'q}$ for
$p \leq p'$ and $s' = s + d/p' -d/p$, one can see that
$f_A$ belongs to $B_{s'2q}$ with $s'=s+d/2-d/p$, and
so by definition, with $\{\epsilon_j\} \in \ell_q$,
$$\|f_A -P_jf_A\|_2 \leq \epsilon_j 2^{-j(s+d/2-d/p)}.$$


If $p > 2$, since we consider compactly supported densities,
with $\{\epsilon_j\} \in \ell_q$,
$$\| f_A-   P_j f_A\|_2 \leq \| f_A-   P_j f_A\|_p 
\leq \epsilon_j 2^{-js}.$$

Finally with $s'$ as claimed, equation \eqref{ketoile} yields
again%
$ |K_\star(A,f)- C_j(f_A)^{1\over2}| \leq  K 2^{-js'}$. 
\endproof


\titre Computation of the estimator $\hat C_j$

The estimator is computable by means of any Daubechies wavelet, including the Haar wavelet.

For a regular wavelet $(D2N, N> 1)$,
it is known how to compute the values $\varphi_{jk}(x)$ (and any derivative) at dyadic rational numbers (Nguyen and Strang, 1996)\nocite{wavelets_and_filter_banks};
this is the approach we have adopted in this paper.

Alternatively, using the customary filtering scheme, one can compute the Haar projection at high $j$ and
use a discrete wavelet transform (DWT) by a $D2N$ to synthetize the coefficients at a lower, more
appropriate resolution before computing the contrast. This avoids the need to precompute any
value at dyadics, because the Haar projection is like a histogram, but adds 
the time of the DWT. 

While this second approach usually gives full satisfaction in density estimation,
in the ICA context, without special care, it can lead to an inflation of computational resources,
or a possibly inoperative contrast at minimization stage. Indeed, 
for the Haar contrast to show any variation
in response to a small perturbation, $j$ must be very high regardless of what would be required by 
the signal regularity and the number of observations; whereas for a D4 and above, 
we just need to set high the precision of dyadic rational
approximation, which present no inconvenience and can be viewed as a memory friendly refined binning inside the binning in $j$.  

We have then chosen the approach with dyadics for simplicity at the minimization stage and possibly
more accurate solutions.

Also for simplicity, in all simulations that follow we have adopted the convention
that the whole signal is contained in the hypercube $[0,1]^{\otimes d}$, after possible rescaling.
For the compactly supported Daubechies wavelets (Daubechies, 1992\nocite{daubechies}), $D2N, N=1,2,\ldots$,
whose support is $[0,2N-1]$,
 the maximum number of $k$ 
 intersecting with an observation lying in the cube is $(2^{j} + 2N-2)^d$.

Note that relocation in the unit hypercube is not a requirement, but otherwise
a sparse array implementation should be used for efficiency.

\subsection Sensitivity of the wavelet contrast

In this section, we compare  independent
and mixed D2 to D8 contrasts on a uniform whitened signal, 
in dimension 2 with $100000$ observations, and 
in dimension 4 with $50000$ observations.
According to proposition \ref{linear-choice}, for $s=+\infty$ the best choice
is $j=0$, to be interpreted as the smallest of technically working $j$, \ie 
satisfying  $2^j > 2N-1$, to ensure that the wavelet support is mostly contained in the observation support. 

For $j=0$, there is only one cell in the cube and the contrast is unable to detect any mixing effect:
for Haar it is identically zero, and for the others D2N it is a constant (quasi for round-off errors)
because we use circular shifting if the wavelet passes an end of the observation support.
At small $j$ such that $2\leq 2^j\leq 2N-1$, D2N wavelets behave more or less like the Haar wavelet, except they are more responsive
to a small perturbation.
We use the Amari distance as defined in Amari (1996) rescaled from 0 to 100\nocite{amaridist}.

In this example, we have deliberately chosen an orthogonal matrix producing 
a small Amari error (less than 1 on a scale from 0 to 100),
pushing the contrast to the limits.
\def\mystrut{\vrule width 0pt height 6pt depth 3pt}
$$
\eqalign{\vbox{\smalltt
\offinterlineskip
\hrule\halign {&\vrule#&\mystrut\quad\hfil#\quad\cr
depth4pt&\omit&&\omit&&\omit&&\omit&\cr
&j&&D2  indep&&D2 mixed&&cpu&\cr
height2.5pt&\omit&&\omit&&\omit&&\omit&\cr
\noalign{\hrule}
height2pt&\omit&&\omit&&\omit&&\omit&\cr
&0 &&0.000E+00&&0.000E+00&&0.12&\cr
&1 &&0.184E-06&&0.102E-10&&0.06&\cr
&2 &&0.872E-04&&0.199E-04&&0.06&\cr
&3 &&0.585E-03&&0.294E-03&&0.06&\cr
&4 &&0.245E-02&&0.285E-02&&0.06&\cr
&5*&&0.926E-02&&0.110E-01&&0.07&\cr
&6 &&0.395E-01&&0.387E-01&&0.07&\cr
&7 &&0.162E+00&&0.162E+00&&0.07&\cr
&8 &&0.651E+00&&0.661E+00&&0.08&\cr
&9 &&0.262E+01&&0.262E+01&&0.12&\cr
&10&&0.105E+02&&0.105E+02&&0.23&\cr
&11&&0.419E+02&&0.419E+02&&0.69&\cr
&12&&0.168E+03&&0.168E+03&&2.48&\cr
height2pt&\omit&&\omit&&\omit&&\omit&\cr
}
\hrule
}
\kern 7pt
\vbox{\smalltt
\offinterlineskip
\hrule\halign {&\vrule#&\mystrut\quad\hfil#\quad\cr
depth4pt&\omit&&\omit&&\omit&&\omit&\cr
&j&&D4  indep&&D4 mixed&&cpu&\cr
height2.5pt&\omit&&\omit&&\omit&&\omit&\cr
\noalign{\hrule}
height2pt&\omit&&\omit&&\omit&&\omit&\cr
&0 &&0.250E+00&&0.250E+00&&0.21&\cr
&1*&&0.239E+00&&0.522E+00&&0.17&\cr
&2 &&0.198E-04&&0.209E-04&&0.17&\cr
&3 &&0.127E-03&&0.159E-03&&0.17&\cr
&4 &&0.635E-03&&0.714E-03&&0.17&\cr
&5 &&0.235E-02&&0.282E-02&&0.17&\cr
&6 &&0.988E-02&&0.105E-01&&0.17&\cr
&7 &&0.405E-01&&0.419E-01&&0.17&\cr
&8 &&0.163E+00&&0.165E+00&&0.21&\cr
&9 &&0.653E+00&&0.653E+00&&0.26&\cr
&10&&0.261E+01&&0.262E+01&&0.39&\cr
&11&&0.104E+02&&0.105E+02&&0.87&\cr
&12&&0.419E+02&&0.420E+02&&2.67&\cr
height2pt&\omit&&\omit&&\omit&&\omit&\cr
}\hrule
}
}$$
\centerline{\smallsf Table 1a. Wavelet contrast values for a D2 and a D4 on a uniform density in dimension 2 
under a half degree rotation }
\centerline{\smallsf Amari error $\approx.8$, nobs=100000,  L=10, }

$$\eqalign{
\vbox{\smalltt
\offinterlineskip
\hrule\halign {&\vrule#&\mystrut\quad\hfil#\quad\cr
depth4pt&\omit&&\omit&&\omit&&\omit&\cr
&j&&D6 indep&&D6 mixed&&cpu&\cr
height2.5pt&\omit&&\omit&&\omit&&\omit&\cr
\noalign{\hrule}
height2pt&\omit&&\omit&&\omit&&\omit&\cr
&0 &&0.304E+00&&0.304E+00&&0.37&\cr
&1 &&0.304E+00&&0.305E+00&&0.37&\cr
&2*&&0.215E+00&&0.666E+00&&0.37&\cr
&3 &&0.132E-03&&0.188E-03&&0.36&\cr
&4 &&0.641E-03&&0.717E-03&&0.36&\cr
&5 &&0.295E-02&&0.335E-02&&0.35&\cr
&6 &&0.123E-01&&0.126E-01&&0.37&\cr
&7 &&0.495E-01&&0.518E-01&&0.36&\cr
&8 &&0.198E+00&&0.200E+00&&0.41&\cr
&9 &&0.796E+00&&0.791E+00&&0.49&\cr
&10&&0.319E+01&&0.319E+01&&0.64&\cr
&11&&0.127E+02&&0.128E+02&&1.13&\cr
&12&&0.509E+02&&0.511E+02&&2.97&\cr
height2pt&\omit&&\omit&&\omit&&\omit&\cr
}\hrule
}
\kern 7pt
\vbox{\smalltt
\offinterlineskip
\hrule\halign {&\vrule#&\mystrut\quad\hfil#\quad\cr
depth4pt&\omit&&\omit&&\omit&&\omit&\cr
&j&&D8 indep&&D8 mixed&&cpu&\cr
height2.5pt&\omit&&\omit&&\omit&&\omit&\cr
\noalign{\hrule}
height2pt&\omit&&\omit&&\omit&&\omit&\cr
&0 &&0.966E+00&&0.966E+00&&0.65&\cr
&1 &&0.966E+00&&0.197E+01&&0.64&\cr
&2*&&0.914E+00&&0.333E+01&&0.65&\cr
&3 &&0.446E-03&&0.409E-03&&0.64&\cr
&4 &&0.220E-02&&0.214E-02&&0.64&\cr
&5 &&0.932E-02&&0.104E-01&&0.63&\cr
&6 &&0.388E-01&&0.383E-01&&0.63&\cr
&7 &&0.157E+00&&0.160E+00&&0.64&\cr
&8 &&0.628E+00&&0.630E+00&&0.71&\cr
&9 &&0.253E+01&&0.252E+01&&0.84&\cr
&10&&0.101E+02&&0.101E+02&&1.03&\cr
&11&&0.405E+02&&0.406E+02&&1.53&\cr
&12&&0.162E+03&&0.162E+03&&3.37&\cr
height2pt&\omit&&\omit&&\omit&&\omit&\cr
}\hrule
}\cr
}
$$
\centerline{\smallsf Table 1b. Wavelet contrast values for a D6 and a D8 on a uniform density in dimension 2 
under a half degree rotation }
\centerline{\smallsf Amari error $\approx.8$, nobs=100000,  L=10, }

First, the Haar contrast is out of touch; at low resolution the mixing passes unnoticed 
because the observations stay in their original bins, 
and at high resolution, as for the other wavelets,
any detection becomes impossible because the ratio $2^{jd}/n$ gets too big, and clearly wanders from the optimal rule of Prop. 
\ref{linear-choice}. 

Had we chosen a mixing with bigger Amari error, say 10, the Haar contrast would have worked  at many more resolutions (this can be checked using the program {\tt icalette1}); still, the Haar contrast is 
less likely to reach small Amari errors in a minimization
process.

For wavelets D4 and above, the contrast is able to capture the mixing
effect especially at low resolution 
(resolution with largest relative increase marked) and up to $j=8$. Also, the wavelet support technical constraint
is apparent between D4 and D6 or D8.
 
Finally we observe that the difference in computing time between
Haar and a D8 is not significative in small dimension;
it gets important starting from dimension 4 (Table 2). Note that the relatively
longer cpu time for  $2^j < 2N-1$ is caused by the need to 
compute a circular shift for practically all points instead of only at borders.

$$\eqalign{
\vbox{\smalltt
\offinterlineskip\hrule\halign {&\vrule#&\mystrut\quad\hfil#\quad\cr
depth4pt&\omit&&\omit&&\omit&&\omit&\cr
&j&&D2 indep&&D2 mixed&&cpu&\cr
height2.5pt&\omit&&\omit&&\omit&&\omit&\cr
\noalign{\hrule}
height2pt&\omit&&\omit&&\omit&&\omit&\cr
&0&&0.000E+00&&0.000E+00&&0.08&\cr
&1&&0.100E-03&&0.155E-06&&0.05&\cr
&2&&0.411E-02&&0.221E-02&&0.05&\cr
&3&&0.831E-01&&0.684E-01&&0.05&\cr
&4&&0.132E+01&&0.129E+01&&0.08&\cr
&5&&0.210E+02&&0.210E+02&&0.29&\cr
&6&&0.336E+03&&0.335E+03&&3.62&\cr
height2pt&\omit&&\omit&&\omit&&\omit&\cr
}\hrule
}
\kern 7pt\vbox{\smalltt

\offinterlineskip\hrule\halign {&\vrule#&\mystrut\quad\hfil#\quad\cr
depth4pt&\omit&&\omit&&\omit&&\omit&\cr
&j&&D4 indep&&D4 mixed&&cpu&\cr
height2.5pt&\omit&&\omit&&\omit&&\omit&\cr
\noalign{\hrule}
height2pt&\omit&&\omit&&\omit&&\omit&\cr
&0&&0.625E-01&&0.625E-01&&0.85&\cr
&1&&0.624E-01&&0.304E+00&&0.83&\cr
&2&&0.283E-03&&0.331E-03&&0.82&\cr
&3&&0.503E-02&&0.453E-02&&0.83&\cr
&4&&0.818E-01&&0.824E-01&&0.92&\cr
&5&&0.130E+01&&0.133E+01&&1.30&\cr
&6&&0.211E+02&&0.211E+02&&4.68&\cr
height2pt&\omit&&\omit&&\omit&&\omit&\cr
}\hrule
}
\cr
\vbox{\smalltt
\offinterlineskip\hrule\halign {&\vrule#&\mystrut\quad\hfil#\quad\cr
depth4pt&\omit&&\omit&&\omit&&\omit&\cr
&j&&D6 indep&&D6 mixed&&cpu&\cr
height2.5pt&\omit&&\omit&&\omit&&\omit&\cr
\noalign{\hrule}
height2pt&\omit&&\omit&&\omit&&\omit&\cr
&0&&0.926E-01&&0.926E-01&&6.03&\cr
&1&&0.927E-01&&0.929E-01&&6.01&\cr
&2&&0.884E-01&&0.825E+00&&6.01&\cr
&3&&0.725E-02&&0.744E-02&&6.07&\cr
&4&&0.122E+00&&0.117E+00&&6.40&\cr
&5&&0.193E+01&&0.195E+01&&7.51&\cr
&6&&0.311E+02&&0.311E+02&&11.0&\cr
height2pt&\omit&&\omit&&\omit&&\omit&\cr
}\hrule
}
\kern 7pt\vbox{\smalltt
\offinterlineskip\hrule\halign {&\vrule#&\mystrut\quad\hfil#\quad\cr
depth4pt&\omit&&\omit&&\omit&&\omit&\cr
&j&&D8 indep&&D8 mixed&&cpu&\cr
height2.5pt&\omit&&\omit&&\omit&&\omit&\cr
\noalign{\hrule}
height2pt&\omit&&\omit&&\omit&&\omit&\cr
&0&&0.934E+00&&0.934E+00&&22.8&\cr
&1&&0.934E+00&&0.364E+01&&22.8&\cr
&2&&0.937E+00&&0.111E+02&&22.8&\cr
&3&&0.751E-01&&0.751E-01&&22.9&\cr
&4&&0.124E+01&&0.117E+01&&24.1&\cr
&5&&0.196E+02&&0.196E+02&&27.0&\cr
&6&&0.313E+03&&0.313E+03&&30.8&\cr
height2pt&\omit&&\omit&&\omit&&\omit&\cr
}\hrule
}
}$$
\centerline{\smallsf  Table 2. Wavelet contrast values on a uniform density, dim=4 , nobs=50000,  L=10, 
 Amari error $\approx.5$}

Computation uses double precision, but single precision works just as well.
There is no guard against inaccurate sums that occur about 10\% of the time
for D4 and above, because it does not prevent a minimum contrast from  
detecting independence. Dyadic approximation parameter $L$ is set 
at octave 10, about three exact decimals, and shows enough. 
Cpu times, in seconds, correspond to the total of the projection time
on $V_j^{d}$ and on the $d$ $V_j$, added to the wavelet contrast computation time;
machine used for simulations is a G4 1,5Mhz, with 1Go ram;
programs are written in fortran and compiled with IBM xlf 
(program {\tt icalette1} to be found in
Appendix).

\subsection Contrast complexity\par

By complexity we mean the length of do-loops.

The projection of $n$ observations on the tensorial space $V_j^d$ and the $d$ margins for a Db(2N)
has complexity $O(n(2N-1)^d)$. 
This is $O(n)$ for a Haar wavelet (2N=2) which boils down to making a histogram. 
The projection complexity is almost independent of $j$ except for memory allocation.
Once the projection at level $j$ is known, the contrast is computed in $O(2^{jd})$.

On the other hand, the complexity to apply one discrete wavelet transform at level $j$ has
complexity $O(2^{jd}(2N-1)^d)$. So we see that 
the filtering approach consisting in taking the Haar projection for a high $j_1$ (typically $2^{j_1d}\approx {n\over \log n}$)
and filter down to a lower $j_0$,
as a shortcut to direct D2N moment approximation at level $j_0$,
is definitely a shortcut; except that the Haar wavelet carries with it 
a lack of sensitivity to small 
perturbations, which is a problem for empirical gradient evaluation
or the detection of a small departure from independence.

For comparison, the Hilbert-Schmidt independence criterion 
is theoretically computed
in $O(n^2)$ (Gretton et al. 2004 section 3), and the Kernel ICA criterion
is theoretically 
computed in $O(n^3d^3)$. In both cases, 
using incomplete Choleski decomposition and 
low-rank approximation of the Gram matrix,
the complexity is brought down in practice to
$O(nd^2)$ for HSIC and  $O(n^2\log n)$ 
for Kernel ICA(Bach and Jordan 2002 p.19).

\titre Contrast minimization 

The natural way to minimize the ICA contrast as a function of a demixing matrix $W$, is
to whiten the signal and then carry out a steepest descent algorithm 
given the constraint $^tWW=I_{d}$, corresponding to $W$ lying on the 
the Stiefel manifold $S(d,d) = O(d)$.
In the ICA context,  we can restrict to $SO(d)\subset O(d)$ 
thus ignoring reflections that are not relevant.  

Needed material for minimization on the Stiefel manifold can be found in the paper of
Arias et al. (1998)\nocite{geometry}. 
Another very close method uses 
the Lie group structure of $SO(d)$ and the corresponding Lie algebra
$so(d)$ mapped together by the matrix logarithm and exponential (Plumbley, 2004)\nocite{plumbley}.
For convenience we reproduce here the algorithm in question,
which is equivalent to a line search in the steepest descent direction in $so(d)$:
\unorderedlist
\li start at $O\in so(d)$, equivalent to $I\in SO(d)$;
\li move about in $so(d)$ from $0$ to $-\eta \nabla_BJ$,
 where $\eta \in \R^+$ corresponds to the minimum in direction $\nabla_BJ$ found by a line search algorithm, 
 where $\nabla_BJ=\nabla J\,{}^tW - W\,{}^t\nabla J$ is the gradient of $J$ in $so(d)$, and where $\nabla J$ is the gradient of $J$ in $SO(d)$;   
\li use the matrix exponential to map back into $SO(d)$, giving $R=\exp(-\eta \nabla_BJ)$;
\li calculate $W'=RW \in SO(d)$ and iterate.
\endunorderedlist

We reproduce below some typical runs (program {\tt icalette3}), with a D4 and
$L=10$. Note that on example 2, the contrast cannot be usefully minimized because of a wrong
resolution.

$$\eqalign{
\vbox{\smalltt
\offinterlineskip\hrule\halign {&\vrule#&\mystrut\quad\hfil#\quad\cr
depth4pt&\omit&\omit&\omit&\omit&\omit&\cr
&\multispan5\hfil d=3, j=3, n=30000 uniform\hfil&\cr
depth3pt&\omit&\omit&\omit&\omit&\omit&\cr
\noalign{\hrule}
depth 3pt&\omit&&\omit&&\omit&\cr
&it&&contrast&&amari&\cr
\noalign{\hrule}
height2pt&\omit&&\omit&&\omit&\cr
&0&&0.127722&&65.842&\cr
&1&&0.029765&&15.784&\cr
&2&&0.002600&&2.129&\cr
&3&&0.001939&&0.288&\cr
&4&&-&&-&\cr
&5&&-&&-&\cr
height2pt&\omit&&\omit&&\omit&\cr
\noalign{\hrule}}
}
\kern 7pt
\vbox{\smalltt
\offinterlineskip\hrule\halign {&\vrule#&\mystrut\quad\hfil#\quad\cr
depth4pt&\omit&\omit&\omit&\omit&\omit&\cr
&\multispan5\hfil d=3, j=5, n=30000 uniform\hfil&\cr
depth3pt&\omit&\omit&\omit&\omit&\omit&\cr
\noalign{\hrule}
depth 3pt&\omit&&\omit&&\omit&\cr
&it&&contrast&&amari&\cr
\noalign{\hrule}
height2pt&\omit&&\omit&&\omit&\cr
&0&&0.321970&&65.842&\cr
&1&&0.321948&&65.845&\cr
&2&&0.321722&&65.999&\cr
&3&&0.321721&&65.999&\cr
&4&&-&&-&\cr
&5&&-&&-&\cr
height2pt&\omit&&\omit&&\omit&\cr\noalign{\hrule}}
}
\kern 7pt
\vbox{\smalltt
\offinterlineskip\hrule\halign {&\vrule#&\mystrut\quad\hfil#\quad\cr
depth4pt&\omit&\omit&\omit&\omit&\omit&\cr
&\multispan5\hfil d=3, j=3, n=10000 uniform\hfil&\cr
depth3pt&\omit&\omit&\omit&\omit&\omit&\cr
\noalign{\hrule}
depth 3pt&\omit&&\omit&&\omit&\cr
&it&&contrast&&amari&\cr
\noalign{\hrule}
height2pt&\omit&&\omit&&\omit&\cr
&0&&0.092920&&42.108&\cr
&1&&0.035336&&14.428&\cr
&2&&0.007458&&3.392&\cr
&3&&0.006345&&1.684&\cr
&4&&0.006122&&1.109&\cr
&5&&0.006008&&0.675&\cr
height2pt&\omit&&\omit&&\omit&\cr
\noalign{\hrule}}
}
\cr
\vbox{\smalltt
\offinterlineskip\hrule\halign {&\vrule#&\mystrut\quad\hfil#\quad\cr
depth4pt&\omit&\omit&\omit&\omit&\omit&\cr
&\multispan5\hfil d=4, j=2, n=10000 uniform\hfil&\cr
depth3pt&\omit&\omit&\omit&\omit&\omit&\cr
\noalign{\hrule}
depth 3pt&\omit&&\omit&&\omit&\cr
&it&&contrast&&amari&\cr
\noalign{\hrule}
height2pt&\omit&&\omit&&\omit&\cr
&0&&0.025193&&22.170&\cr
&1&&0.010792&&9.808&\cr
&2&&0.003557&&4.672&\cr
&3&&0.001272&&1.167&\cr
&4&&0.001033&&0.502&\cr
&5&&0.000999&&0.778&\cr
height2pt&\omit&&\omit&&\omit&\cr
\noalign{\hrule}}
}
\kern 7pt
\vbox{\smalltt
\offinterlineskip\hrule\halign {&\vrule#&\mystrut\quad\hfil#\quad\cr
depth4pt&\omit&\omit&\omit&\omit&\omit&\cr
&\multispan5\hfil d=3, j=4, n=30000 expone.\hfil&\cr
depth3pt&\omit&\omit&\omit&\omit&\omit&\cr
\noalign{\hrule}
depth 3pt&\omit&&\omit&&\omit&\cr
&it&&contrast&&amari&\cr
\noalign{\hrule}
height2pt&\omit&&\omit&&\omit&\cr
&0&&8.609670&&52.973&\cr
&1&&5.101633&&48.744&\cr
&2&&0.778619&&16.043&\cr
&3&&0.017585&&3.691&\cr
&4&&0.008027&&2.262&\cr
&5&&0.006306&&1.542&\cr
height2pt&\omit&&\omit&&\omit&\cr\noalign{\hrule}}
}
\kern 7pt
\vbox{\smalltt
\offinterlineskip\hrule\halign {&\vrule#&\mystrut\quad\hfil#\quad\cr
depth4pt&\omit&\omit&\omit&\omit&\omit&\cr
&\multispan5\hfil d=3, j=3, n=10000 semici.\hfil&\cr
depth3pt&\omit&\omit&\omit&\omit&\omit&\cr
\noalign{\hrule}
depth 3pt&\omit&&\omit&&\omit&\cr
&it&&contrast&&amari&\cr
\noalign{\hrule}
height2pt&\omit&&\omit&&\omit&\cr
&0&&0.041392&&35.080&\cr
&1&&0.029563&&22.189&\cr
&2&&0.007775&&5.601&\cr
&3&&0.006055&&3.058&\cr
&4&&0.005387&&2.261&\cr
&5&&0.005355&&1.541&\cr
height2pt&\omit&&\omit&&\omit&\cr\noalign{\hrule}}
}
}$$
\centerline{\smallsf Table 3. Minimization examples at various j, d and n with D4 and L=10 }%

In our simulations, $\nabla J$ is computed by first differences; 
in doing so we cannot keep perturbed $W$ orthogonality, 
and we actually compute a plain gradient in $\R^{dd}$. 

Again,  a Haar contrast empirical gradient 
is tricky to obtain, since a small perturbation in $W$ will 
likely result in an unchanged histogram at small $j$, whereas 
with D4 and above contrasts, response to perturbation 
is practically automatic and is anyway adjustable by the dyadic 
approximation parameter $L$.\nocite{EXPOKIT}

Below is the average of 100 runs in dimension 2 with 10000 observations, D4,
$j=3$ and $L=10$ for different densities;  {\tt start} columns indicate
Amari distance (on the scale 0 to 100) and  wavelet contrast on entry; {\tt it}
column is the average number of iterations. Note that for some densities after whitening we are already close to the minimum, but the contrast still
detects a departure from independence; the routine exits on entry
if the contrast or the gradient are too small, and this practically always correspond to an Amari distance less than 1 in our simulations.
$$\vbox{\smalltt
\offinterlineskip\hrule\halign {&\vrule#&\mystrut\quad\hfil#\quad\cr
depth 4pt&\omit&&\omit&&\omit&&\omit&&\omit&&\omit&\cr
&density&&Amari start&&Amari end&&cont. start&&cont. end&&it.&\cr
height 2pt&\omit&&\omit&&\omit&&\omit&&\omit&&\omit&\cr
\noalign{\hrule}
height 2pt&\omit&&\omit&&\omit&&\omit&&\omit&&\omit&\cr
&uniform&&53.193&&0.612&&0.509E-01&&0.104E-02&&1.7&\cr
&exponential&&32.374&&0.583&&0.616E-01&&0.150E-03&&1.4&\cr
&Student&&2.078&&1.189&&0.534E-04&&0.188E-04&&0.1&\cr
&semi-circ&&51.401&&2.760&&0.222E-01&&0.165E-02&&1.8&\cr
&Pareto&&4.123&&0.934&&0.716E-03&&0.415E-05&&0.3&\cr
&triangular&&46.033&&7.333&&0.412E-02&&0.109E-02&&1.6&\cr
&normal&&45.610&&45.755&&0.748E-03&&0.408E-03&&1.4&\cr
&Cauchy&&1.085&&0.120&&0.261E-04&&0.596E-06&&0.1&\cr
height 2pt&\omit&&\omit&&\omit&&\omit&&\omit&&\omit&\cr
\noalign{\hrule}
}}  
$$
\centerline{\smallsf  Table 4. Average results of 100 runs in dimension 2, j=3
with a D4 at L=10}

These first results are comparable with the performance of existing ICA algorithms,
as presented for instance in the paper of Jordan and Bach (2002) p.30
(average Amari error between 3 and 10 for 2 
sources and 1000 observations)
or Gretton et al (2004) table 2 (average Amari error between 2 and 6 for 2 
sources and 1000 observations).

Finally we give other runs on the example of the uniform density
at resolution $j=2$ under different parameters settings, and relatively
fewer number of observations.

$$\vbox{\smalltt
\offinterlineskip\hrule\halign {&\vrule#&\mystrut\quad\hfil#\quad\cr
depth 4pt&\omit&&\omit&&\omit&&\omit&&\omit&&\omit&&\omit&&\omit&\cr
&obs.&&dim&&L&&Amari start&&Amari end&&cont. start&&cont. end&&it.&\cr
height 2pt&\omit&&\omit&&\omit&&\omit&&\omit&&\omit&&\omit&&\omit&\cr
\noalign{\hrule}
height 2pt&\omit&&\omit&&\omit&&\omit&&\omit&&\omit&&\omit&&\omit&\cr
&250  &&2 &&10 &&47.387  &&38.919 &&0.279E-01 &&0.193E-01   &&2.4&\cr
&250  &&2 &&13 &&47.387  &&32.470 &&0.279E-01 &&0.170E-01   &&2.2&\cr
&250  &&2 &&16 &&47.387  &&17.915 &&0.279E-01 &&0.603E-02   &&2.3&\cr
&250  &&2 &&19 &&47.387  &&19.049 &&0.279E-01 &&0.598E-02   &&2.6&\cr
&500  &&2 &&10 &&51.097  &&20.700 &&0.246E-01 &&0.106E-01   &&2.1&\cr
&500  &&2 &&13 &&51.097  && 6.644 &&0.246E-01 &&0.398E-02   &&2.2&\cr
&500  &&2 &&16 &&51.097  &&21.063 &&0.246E-01 &&0.109E-01   &&2.1&\cr
&500  &&2 &&19 &&51.097  &&14.734 &&0.246E-01 &&0.839E-02   &&2.4&\cr
&1000 &&2 &&10 &&41.064  && 3.533 &&0.167E-01 &&0.186E-02   &&2.3&\cr
&1000 &&2 &&13 &&41.064  && 3.071 &&0.167E-01 &&0.190E-02   &&2.1&\cr
&1000 &&2 &&16 &&41.064  && 3.518 &&0.167E-01 &&0.194E-02   &&1.9&\cr
&1000 &&3 &&16 &&49.607  &&15.082 &&0.405E-01 &&0.127E-01   &&4.8&\cr
&5000 &&3 &&10 &&49.575   &&5.405 &&0.390E-01 &&0.399E-02   &&4.5&\cr
&5000 &&3 &&16 &&49.575   &&1.668 &&0.390E-01 &&0.960E-03   &&4.7&\cr
&5000 &&4 &&10 &&43.004  &&17.036 &&0.561E-01 &&0.190E-01   &&4.4&\cr
&5000 &&5 &&10 &&38.400  &&29.679 &&0.800E-01 &&0.559E-01   &&4.1&\cr
&5000 &&5 &&16 &&38.400  && 4.233 &&0.798E-01 &&0.700E-02   &&5.0&\cr
&5000 &&6 &&16 &&42.529  &&10.841 &&0.114E+00 &&0.278E-01   &&4.9&\cr
&5000 &&7 &&16 &&41.128  &&15.761 &&0.188E+00 &&0.573E-01   &&5.0&\cr
&5000 &&8 &&16 &&39.883  &&14.137 &&0.286E+00 &&0.743E-01   &&5.0&\cr
height 2pt&\omit&&\omit&&\omit&&\omit&&\omit&&\omit&&\omit&&\omit&\cr
\noalign{\hrule}
}}  
$$
\centerline{\smallsf  Table 5. Average results of 10 runs, j=2,
with a D4, truncated at 5 iterations. }

One can see that raising the dyadic approximation parameter $L$ tends
to improve the minimization when the number of observations is "low" relatively
to the number or cells $2^{jd}$, but that 500 observations in dimension 2 seems
to be a minimum in the current state of the program. In higher dimensions,
a higher number of observations is required, and in dimension 6 and above, 5000
is not enough at L=16.

\subsection A visual example in dimension 2\par

In dimension 2, we are exempted from any added complication brought by a gradient
descent and Stiefel minimization.
After whitening, the inverse of A is an orthogonal matrix, whose membership can be restricted
to $SO(2)$, ignoring reflections. 
So there is only one parameter
$\theta$ to find to achieve reverse mixing. Since permutations of axes are also void operations
in ICA, angles in the range 0 to $\pi/2$ are enough to find out the minimum $W_0$ which, right
multiplied by N, will recover the ICA inverse of A. And A can be set to the identity matrix,
because what changes when A is not the identity, but any invertible matrix, is completely
contained in N.

Figures below show the wavelet contrast in W and the amari distance $d(A,WN)$ (where
N is the matrix computed after whitening),  functions of the rotation angle of the
matrix $W$ restricted to one period, $[ 0, \pi/2 ]$. The minimums are not necessarily at a zero
angle, for precisely, mere whitening leaves the signal in a 
random rotated position (to reproduce the
following results run the program {\tt icalette2}).

We see that, provided Amari error and wavelet contrast have coinciding
minima, any line search algorithm will find the angle to reverse the mixing effect. We see also in Fig.2 that the Haar wavelet contrast 
is perfectly suitable to detect independence, so that minimization methods not gradient based might work very well in this case.

On the example of the uniform density (Fig.3) we have an illustration
of  a non smooth contrast variation typical of a too high resolution $j$ given
regularity and number of observations. 
$$
\vbox{
\setbox0=\hbox{\epsfxsize=5cm \epsfbox{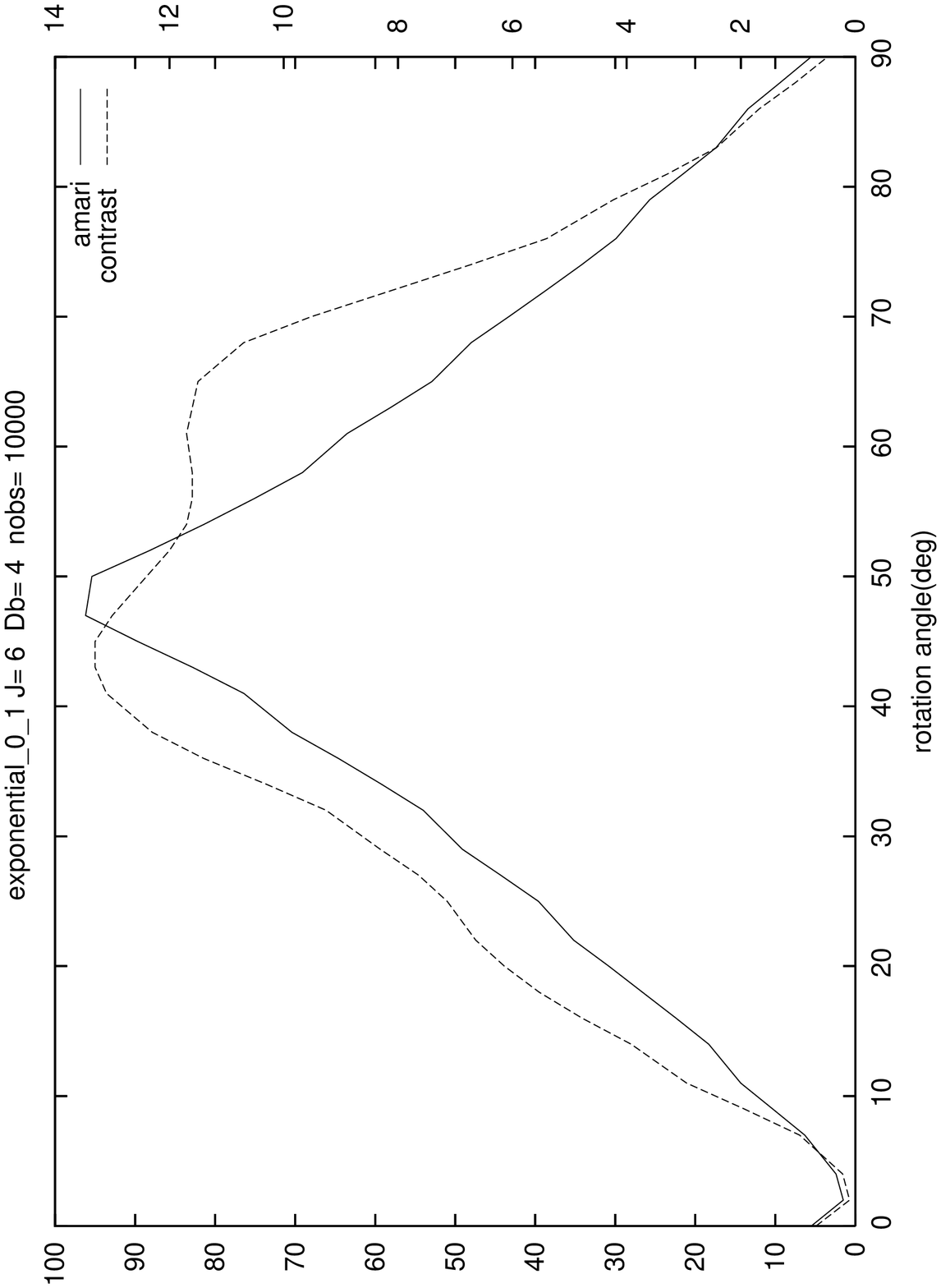}}
\rotr0
\box0\hbox to 7cm{\hfil\smallsf Fig.1. Exponential, D4, j=6, n=10000\hfil}
}
\vbox{
\setbox0=\hbox{\epsfxsize=5cm \epsfbox{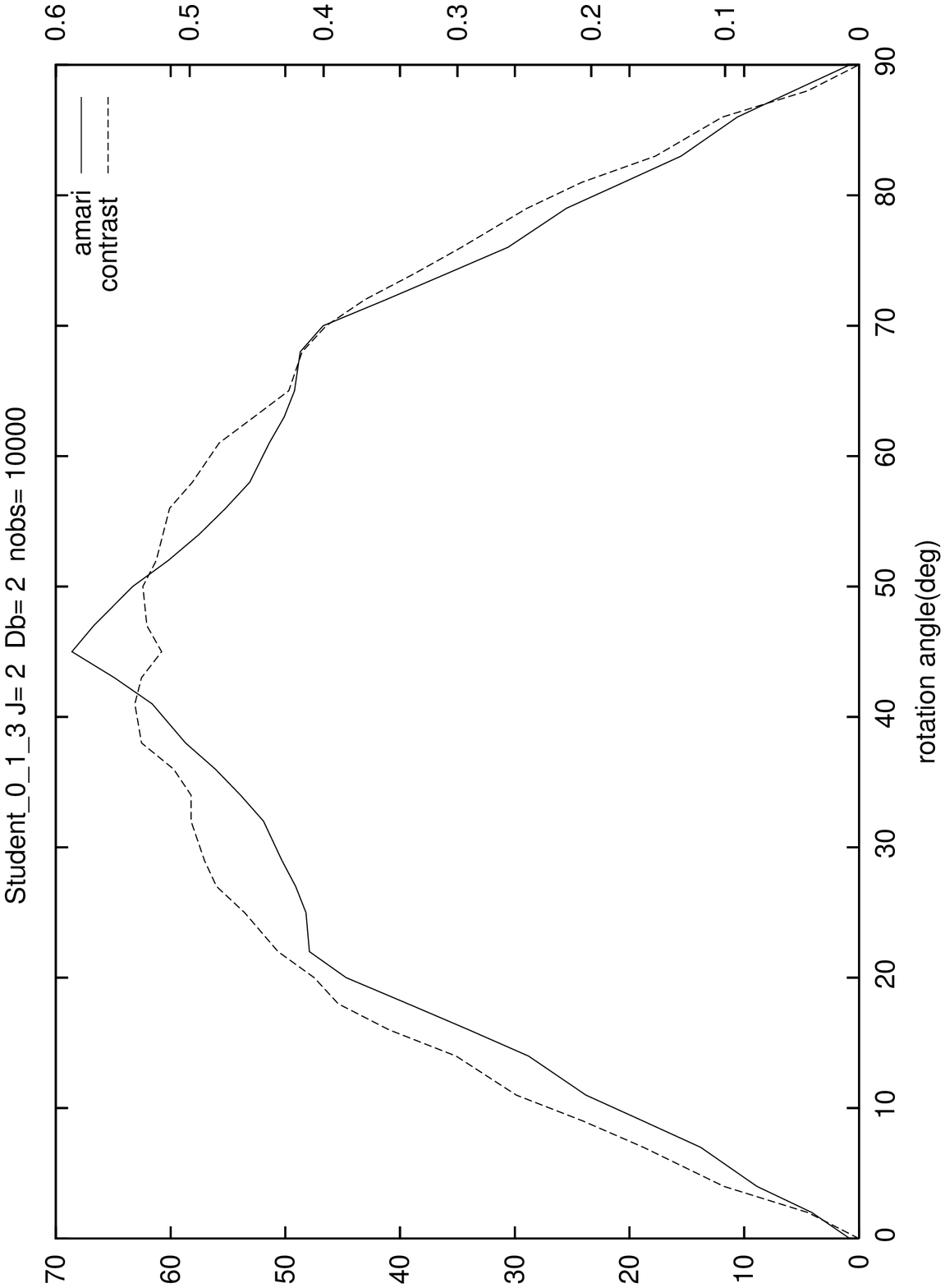}}
\rotr0
\box0\hbox to 7cm{\hfil\smallsf Fig.2. Student, D2, j=2, n=10000\hfil}
}
$$

$$
\vbox{
\setbox0=\hbox{\epsfxsize=5cm \epsfbox{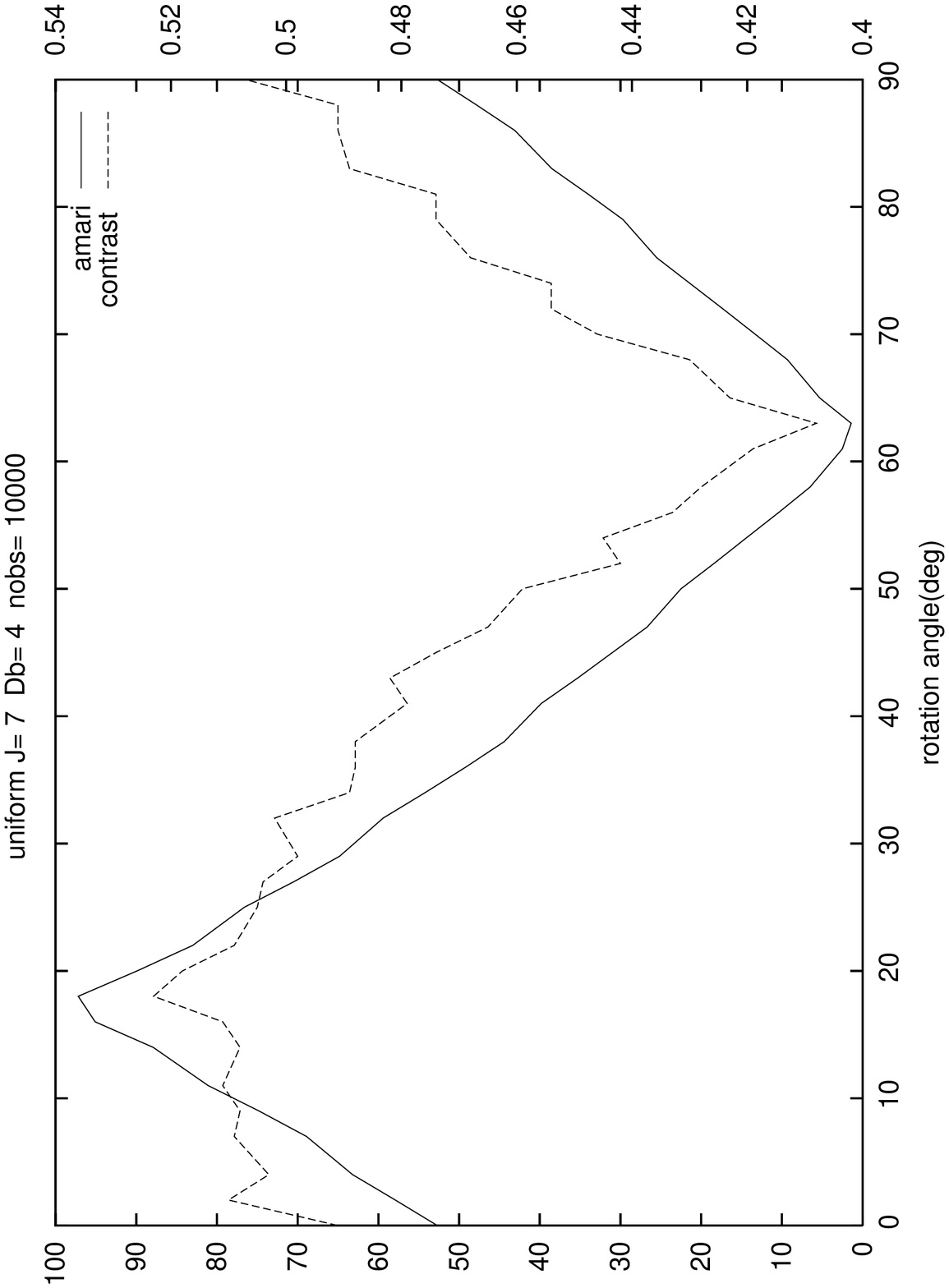}}
\rotr0
\box0\hbox to 7cm{\hfil\smallsf Fig.3. Uniform, D4, j=7, n=10000\hfil}}
\vbox{
\setbox0=\hbox{\epsfxsize=5cm \epsfbox{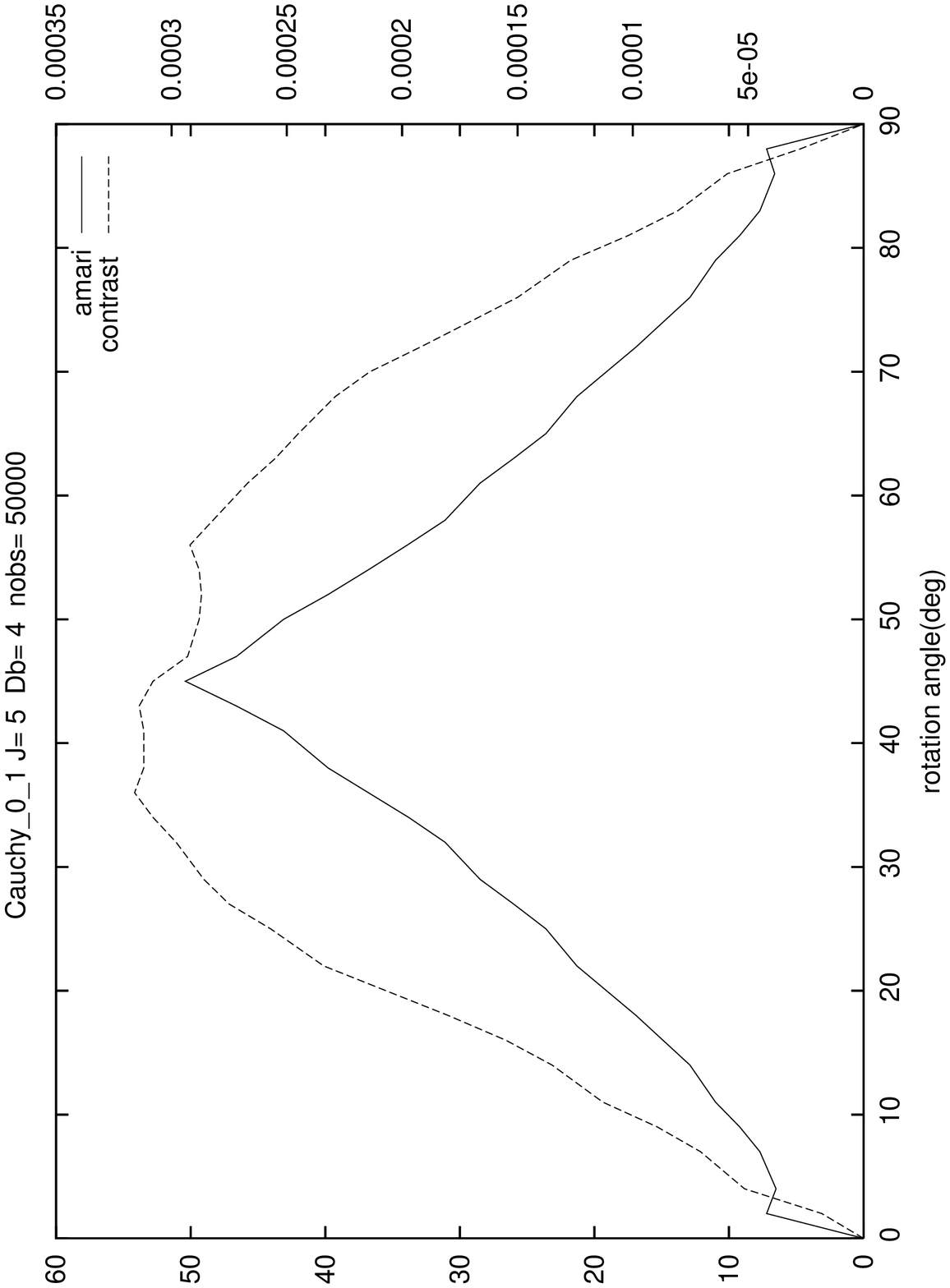}}
\rotr0
\box0\hbox to 7cm{\hfil\smallsf Fig.4. Cauchy, D4, j=5, n=50000\hfil }}
$$

\titre Appendix

\lemma Approximation by the $j$-dependent U-statistic\par
\definexref{UV}{\the\notitre.\the\nolemma}{lemma}
{\it

Let $(\echn(X))$ be an i.i.d. sample of a random variable on $\R^d$
and let 
$\rho,\, \sigma, \, m$ be positive integers verifying $\rho+\sigma=m$.
Let $h_j$ be a (possibly unsymmetric) kernel function on $\R^{md}$
defined as
 $$ h_j(x_1,\ldots,x_m)= \Phi_{jk}(x_1)\ldots \Phi_{jk}(x_\rho)
 \varphi_{jk^{\ell_1}}\circ\pi^{\ell_1}(x_{\rho+1})
 \ldots\varphi_{jk^{\ell_\sigma}}\circ\pi^{\ell_\sigma}(x_{\rho+\sigma}),$$
 with $\pi^\ell$ the canonical projection on component $\ell$
 and $\Phi(x)=\prod_{\ell=1}^d\varphi\circ\pi^\ell(x)$.
  
Consider the associated U-statistic and Von Mises V statistic, 
$$U_{nj}={(n-m)!\over n!}\sum_{i_1\neq\ldots\neq i_m} h_j(X_{i_1},\ldots, X_{i_m}),\quad\qquad
V_{nj}={1\over n^m}\sum_{i_1,\ldots, i_m} h_j(X_{i_1},\ldots,X_{i_m}). 
$$

The following relation holds:
$$E |U_{nj} - V_{nj} |^r = 
\left({2^{j (d\rho+\sigma)}\over n }\right)^{{r\over2}-1} O(n^{-1-{r\over 2}}).\eqdef{ustatr}$$

In corollary, $E |U_{nj} - V_{nj} |^r = O(n^{-r})$ for $r=1,\, 2$ and
$E |U_{nj} - V_{nj} |^r = O(n^{-1-{r\over 2}})$
for $r\geq2$ and $2^{j(d\rho+\sigma)}/n <1$.
}

\proof

The first lines use the proof of the original lemma (see for example Serfling, 1980), with special care for unsymmetric kernels.

Let $W_{nj}$ be the average of all terms $h_j(X_{i_1},\ldots,X_{i_m})$ 
with at least one equality $i_a=i_b$ for $a\neq b$ and $1\leq a, b \leq m$;
there are $n^m-A_n^m$ such terms. 

$\Omega$ denoting the set of $n^m$ unconstrained indexes, one has
the relation
$$\eqalign{
(n^m-A_n^m)(W_{nj}-U_{nj})&=\sum_{\Omega-\{i_1\neq \ldots\neq i_m\}} h(X_{i_1},\ldots,X_{i_m})
- n^m U_{nj} + \sum_{\{i_1\neq \ldots\neq i_m\}} h(X_{i_1},\ldots,X_{i_m})\cr
&=n^m(V_{nj}-U_{nj})
}$$

Hence, using Minkowski inequality and the fact that 
$(n^m-A_n^m)=O(n^{m-1})$ is positive, one obtains,
$$\eqalign{
\Enfa |U_{nj} - V_{nj}|^r
& = n^{-mr}(n^m - A_n^m)^r\,  \Enfa |U_{nj}-W_{nj}|^r\cr
& \leq n^{-mr}(n^m - A_n^m)^r \left(\Enfa |U_{nj}|^r + \Enfa |W_{nj}|^r\right)\cr
& \leq O(n^{-r}) \left(\Enfa |U_{nj}|^r + \Enfa |W_{nj}|^r\right).\cr
}$$
It remains to bound to right parenthesis.

Using Minkowski inequality, one has,%
$$\eqalign{
\Enfa{}^{1\over r} |U_n|^r &\leq 
 [A_n^m]^{-1} \sum_{i_1 \neq\ldots\neq i_m} 
\Enfa{}^{1\over r} |
\Phi_{jk}(X_{i_1})\ldots
\Phi_{jk}(X_{i_\rho})
\varphi_{jk^{\ell_1}}(X_{i_{\rho+1}}^{\ell_{1}})
\ldots\varphi_{jk^{\ell_{i_{\sigma}}}}(X_{i_{m}}^{\ell_{i_{\sigma}}})|^r\cr
&=  \Enfa{}^{1\over r} |\varphi_{jk^{\ell_1}}(X_1^{\ell_1})|^r
\ldots\Enfa{}^{1\over r} |\varphi_{jk^{\ell_\sigma}}
(X_1^{\ell_{i_\sigma}})|^r
\Enfa{}^{\rho\over r} |\Phi_{jk}(X_1)|^{r}
 \cr
}$$ 

Next 
$\Enfa|\varphi_{jk^\ell}(X_1^\ell)|^r = 2^{jr/2}\int |\varphi(2^jx-k^\ell)|^r\fix(f_A,\ell)(x)dx\leq 2^{j({r\over2}-1)} \|\fix(f_A,\ell)\|_\infty \|\varphi\|_r^r$ and by the same means,
$\Enfa|\Phi_{jk}(X_i)|^r\leq 
 2^{jd({r\over2}-1)} \|f_A\|_\infty \|\Phi\|_r^r$

So that,
$$\eqalign{
\Enfa |U_n|^r 
&\leq 2^{j\sigma({r\over2}-1)}
\bigg(\prod_{i=1,\ldots,\sigma}\|\fix(f_A,\ell_{i})\|_\infty\bigg)
\|\varphi\|_r^{r\sigma}\quad
2^{jd\rho({r\over2}-1)}
\|f_A\|_\infty^\rho
\|\Phi\|_r^{r\rho}\cr
&= C 2^{j(\sigma+d\rho)({r\over2}-1)} 
} 
$$

Likewise for $W_n$ one obtains terms of the type
$\Enfa \varphi_{jk^{\ell_1}}(X_1^{\ell_1})\ldots
 \varphi_{jk^{\ell_\kappa}}(X_1^{\ell_\kappa})$ of which $\Enfa\Phi(X_1)$
 is one particular form, which are bounded exactly in the same way,
 and produce the same power of $2^j$.
\endproof

\bigskip\bigskip
I am mostly grateful to my advisor, Dominique Picard, for  many suggestions in
the writing of this paper.

{\bf References}\par

\tenrm
\parskip=10pt

(Achard et al. 2003)
Christian~Jutten S.~Achard, D.T.Pham.
 A quadratic dependence measure for nonlinear blind sources
  separation.
 {\it Proceeding of ICA 2003 Conference}, pages 263--268, 2003.

(Amari, 1996)
A.~Cichocki S.~Amari and H.~Yang.
 A new algorithm for blind signal separation.
 {\it Advances in Neural Information Processing Systems}, 8:757--763,
  1996.

(Arias et al. 1998)
Steven T.~Smith Alan~Edelman, Tomas~Arias.
 {\it The geometry of algorithms with orthogonality constraints}.
 SIAM, 1998.


(Bach \& Jordan, 2002)
M.~I.~Jordan F.~R.~Bach.
 Kernel independent component analysis.
 {\it J. of Machine Learning Research}, 3:1--48, 2002.

(Bell \& Sejnowski, 1995) A.~J.~Bell. T.J.~Sejnowski 
 A non linear information maximization algorithm that performs blind
  separation.
 {\it Advances in neural information processing systems}, 1995.

(Bergh \& L\"ofstr\"om, 1976)
J.~Bergh and J.~L{\"o}str\"om.
 {\it Interpolation spaces}.
 Springer, Berlin, 1976.


(Cardoso, 1999) 
J.F. Cardoso.
 High-order contrasts for independent component analysis.
 {\it Neural computations 11}, pages 157--192, 1999.

(Comon, 1994) 
P.~Comon.
 Independent component analysis, a new concept ?
 {\it Signal processing}, 1994.

(Daubechies, 1992)
Ingrid Daubechies.
 {\it Ten lectures on wavelets}.
 SIAM, 1992.

(Devore \& Lorentz, 1993)
R.~Devore, G.~Lorentz. {\it Constructive approximation}. Springer-Verlag,
1993.

(Donoho et al., 1996)
G.~Kerkyacharian D.L.~Donoho, I.M.~Johnstone and D.~Picard.
 Density estimation by wavelet thresholding.
 {\it Annals of statistics}, 1996.

(Gretton et al. 2003)
Alex~Smola Arthur~Gretton, Ralf~Herbrich.
 The kernel mutual information.
 Technical report, Max Planck Institute for Biological Cybernetics,
  April 2003.

(Gretton et al. 2004)
A.~Gretton,  O.~Bousquet, A.~Smola, B.~Sch\"o{}lkopf.
 Measuring statistical dependence with Hilbert-Schmidt norms.
 Technical report, Max Planck Institute for Biological Cybernetics,
  October 2004.

(H\"ardle et al., 1998)
 Wolfgang H{\"
  a}rdle, G\'erard~Kerkyacharian, Dominique~Picard and 
 Alexander~Tsybakov.
 {\it Wavelets, approximation and statistical applications}.
 Springer, 1998.

(Hyvar\"\i nen et al. 2001)
A.~Hyvar\"\i nen, J.~Karhunen. E.~Oja
{\it Independent component analysis}.
 Inter Wiley Science, 2001.

(Hyvarinen \& Oja, 1997)
A.~Hyvarinen and E.~Oja.
 A fast fixed-point algorithm for independent component analysis.
 {\it Neural computation}, 1997.

(Kerkyacharian \& Picard, 1992)
G{\'e}rard~Kerkyacharian Dominique~Picard.
 Density estimation in {B}esov spaces.
 {\it Statistics and Probability Letters}, 13:15--24, 1992.

(Meyer, 1997)
Yves Meyer.
 {\it Ondelettes et op\'erateurs}.
 Hermann, 1997.


(Nguyen \ Strang, 1996)
Truong~Nguyen Gilbert~Strang.
 {\it Wavelets and filter banks}.
 Wellesley-Cambridge Press, 1996.

(Nikol'ski$\breve{\hbox{\i}}$, 1975) S.M. 
Nikol'ski$\breve{\hbox{\i}}$. Approximation of functions of several variables and imbedding theorems. {\it Springer Verlag, 1975}.

(Peetre, 1975)
Peetre, J.  
New Thoughts on Besov Spaces. Dept. Mathematics, Duke Univ, 1975.

(Plumbley, 2004)
Mark~D. Plumbley.
 Lie group methods for optimization with orthogonality constraints.
 {\it Lecture notes in Computer science}, 3195:1245--1252, 2004.


(Roch, 1995) Jean Roch.
 Le mod\`ele factoriel des cinq grandes dimensions de personnalit\'e :
  les big five.
 Technical report, AFPA, DO/DEM, March 1995.

(Rosenblatt, 1975)
M.~Rosenblatt.
 A quadratic measure of deviation of two-dimensional density estimates
  and a test for independence.
 {\it Annals of Statistics}, 3:1--14, 1975.

(Rosenthal, 1972)
Rosenthal, H. P. 
On the span in lp of sequences of independent random variables. {\it Israel
J. Math}. 8 273–303, 1972.

(Serfling, 1980)
Robert~J. Serfling.
 {\it Approximation theorems of mathematical statistics}.
 Wiley, 1980.

(Sidje, 1998)
R.~B. Sidje.
 {\sc Expokit.} {A} {S}oftware {P}ackage for {C}omputing {M}atrix
  {E}xponentials.
 {\it ACM Trans. Math. Softw.}, 24(1):130--156, 1998.


(Tsybakov \& Samarov,  2004) 
A.~Tsybakov A.~Samarov.
Nonparametric independent component analysis.
 {\it Bernouilli}, 10:565--582, 2004.

(Triebel, 1992)
Triebel, H. Theory of Function Spaces 2. Birkh\"auser, Basel, 1992

{\bf Programs and other runs available at \tt http://www.proba.jussieu.fr/pageperso/barbedor}

\bye